\documentclass{article}
\usepackage[dvipdfmx]{graphicx}
\usepackage[top=30truemm,bottom=30truemm,left=25truemm,right=25truemm]{geometry}
\usepackage{moreverb}
\usepackage{url}
\usepackage{grffile}
\usepackage{lineno}
\usepackage{lipsum}
\usepackage[UKenglish]{isodate}
\usepackage{color,soul}

\usepackage{amsmath,amssymb,amsfonts,mathtools,mathrsfs,bm,url}
\usepackage{psfrag,epsf}
\usepackage{enumerate}
\usepackage{natbib}
\usepackage{tabu}

\newtheorem{theorem}{Theorem}

\newtheorem{corol}{Corollary}
\newtheorem{assumption}{Assumption}
\newtheorem{remark}{Remark}
\newtheorem{proposition}{Proposition}

\title{Parametric and nonparametric probability distribution estimators of sample maximum}
\author{
Taku MORIYAMA\\
School of Data Science, Yokohama City University}
\date{}

\begin{document}
\maketitle

\begin{abstract}
Extreme value theory has constructed asymptotic properties of the sample maximum. This study concerns probability distribution estimation of the sample maximum. The traditional approach is parametric fitting to the limiting distribution---the generalized extreme value distribution; however, the model in non-limiting cases is misspecified to a certain extent. We propose a plug-in type of nonparametric estimator that does not need model specification. Asymptotic properties of the distribution estimator are derived. The simulation study numerically investigates the relative performance in finite-sample cases.

This study assumes that the underlying distribution of the original sample belongs to one of the Hall class, the Weibull class or the bounded class, whose types of the limiting distributions are all different: the Fr\'{e}chet, Gumbel or Weibull. It is proven that the convergence rate of the parametric fitting estimator depends on both the extreme value index and the second-order parameter and gets slower as the extreme value index tends to zero. On the other hand, the rate of the nonparametric estimator is proven to be independent of the extreme value index under certain conditions. The numerical performances of the parametric fitting estimator and the nonparametric estimator are compared, which shows that the nonparametric estimator performs better, especially for the extreme value index close to zero. Finally, we report two real case studies: the Potomac River peak stream flow (cfs) data and the Danish Fire Insurance data.  
\end{abstract}

\section{Introduction}

In many applications, not the mean but extreme values are of interest. In particular, it is important to measure risk accurately by considering not the upper limit but the sample maximum of the appropriate forecast period. In Hydrology and Meteorology, the Potomac River peak stream flow (cfs) data for water in 1895 through 2000 (Oct--Sep) at Point Rocks, Maryland was offered by the US Geological Survey. It is available as Potomac in the \texttt{extRemes} package in the R software environment. The Danish Fire Insurance data were supplied by Mette Rytgaard of Copenhagen Re (see \cite{resnick1997discussion}). \cite{beirlant1996practical} surveyed the Norwegian fire insurance data. In reliability theory, \cite{komukai1994requirements} estimated the life prediction of copper pipes based on maximum pit depths. In corrosion science, \cite{kasai2016predicting} predicted the actual maximum depths of corrosion using detected depths. The sample maximum is sometimes an outlier that is detected by estimating the underlying distribution (\citealt{Mittnik2001distribution}, \citealt{gbenro2020using}).

Let $X_1,X_2,\cdots, X_{n+m}$ be independent and identically distributed random variables with a continuous distribution function $F$. Here, we consider estimating the sample maximum distribution (SMD) $F^m(x)$ of the observations $\{X_{n+1}, \cdots, X_{n+m}\}$ from $\{X_{1}, \cdots, X_{n}\}$. This is useful and provides important information about weather, natural disasters, financial markets, and many other practical phenomena. Estimating the return level needs estimated SMD. We suppose that $n$ is sufficiently large. To investigate asymptotic properties, we regard $m$ as a function of $n$ according to research in extreme value theory i.e. $m=m_n \to \infty$ throughout the paper and suppose $m/n=o(1)$. \cite{novak2009impossibility} proved that there are no consistent estimators of SMD unless $m/n=o(1)$ without any assumptions. We consider the pointwise estimation of SMD at the point $x=x_n \to \infty$ as $n \to \infty$. 

By applying the Fisher-Tippett-Gnedenko theorem, SMD at the point $x$ is approximated by the generalized extreme value distribution (GEV) $G_{\bm{\gamma}}(x)$ for large enough $m$, where 
\begin{align*}
G_{\bm{\gamma}}(x) :=& G_{\gamma}\left(\frac{x -b_m}{a_m}\right) ~~~ \text{for} ~~ 1+{\gamma} \frac{x -{b}_m}{{a}_m} > 0, \\
G_{\gamma}(x) :=&
\begin{dcases}
\exp(-(1+\gamma x)^{-1/\gamma}) ~~~ &\text{for} ~~ 1+\gamma x>0 ~~~\text{and} ~~~ \gamma \in \mathbb{R} \setminus \{0\} \\
\exp(-\exp(-x)) ~~~ &\text{for} ~~ x \in \mathbb{R} ~~~\text{and} ~~~ \gamma =0.
\end{dcases}
\end{align*}

This study investigates asymptotic properties of two different approaches in estimation of SMD. The fitting estimator works only if the Fisher-Tippett-Gnedenko theorem holds (it is called as $F$ belonging to the maximum domain of attraction). The unknown parameters are estimated by maximum likelihood method based on the block maxima. The details are given in the Appendices. 

In order to obtain the explicit forms of asymptotic properties we need to specify the form of $F$, and to make a theoretical comparison it is necessary that $F$ belongs to the maximum domain of attraction. We suppose that $F$ belongs to either one of the following three types of distribution: the Hall class of distributions (\citealt{hall1984best}), which satisfies $\exists(\alpha, \beta, A, B)$ s.t. $\alpha>0, ~~~  \beta\ge2^{-1}, ~~~  A>0, ~~~  B\neq0$ and
\begin{align}
x^{\alpha+\beta} \{1 - F(x) - A x^{-\alpha} (1 + B x^{-\beta} ) \} \to 0~~~ \text{as}  ~~~ x \to \infty,
\end{align}
or the Weibull class of distributions, which satisfies $\exists(\kappa, C)$ s.t. $\kappa>0, ~~~  C>0$ and
\begin{align}
\exp(Cx^{\kappa})\{1 - F(x) - \exp(-Cx^{\kappa})\} \to 0~~~ \text{as}  ~~~ x \to \infty,
\end{align}
or the bounded class of distributions (see e.g. \cite{stupfler2016estimating}), which satisfies $\exists(x^*, \mu, \sigma, D, E)$ s.t. $x^*\in \mathbb{R}, \mu<-2, \sigma\le-2^{-1},  D>0,  E\neq0$ and
\begin{align}
(x^* -x)^{\mu+\sigma}\{1 -F(x) - (x^* -x)^{-\mu} (D + E (x^* -x)^{-\sigma})\} \to 0 ~~~ \text{as}  ~~~ x \uparrow x^*.
\end{align}

\begin{remark}
The Pareto distributions, $t$-distributions, Burr distributions and extreme value distributions with $\gamma>0$ belong to the Hall class. The distributions belonging to the Hall class are heavy-tailed. Methods for checking whether the underlying distribution is heavily tailed are provided in \cite{resnick1997discussion}. The Weibull distributions and logistic distributions belong to the Weibull class. The distributions belonging to the bounded class each have their upper bounded endpoint. Examples include the uniform distributions, beta distributions and inverse Burr distribution with $c>0$ and $\ell>0$ defined as follows:  
\begin{align*}
1 -F(x) =& (1 - \{1+(-x)^{\ell}\}^{-1})^{c} ~~~ \text{for} ~~ x<0 \\
=& (-x)^{c\ell} (1 + O((-x)^{\ell})) ~~~ \text{as}  ~~~ x \uparrow 0
\end{align*}
(see \citealt{kleiber2003statistical}). The inverse Burr distribution is also called as the Dagum distribution (see also \citealt{kleiber2008guide}).
\end{remark}

Let 
$$M_n := m\times
\begin{dcases}
A x^{-\alpha} ~~~ &\text{for the Hall class} \\
\exp(-C x^{\kappa}) ~~~ &\text{for the Weibull class} \\
D (x^* -x)^{-\mu} ~~~ &\text{for the bounded class},
\end{dcases}$$ 
where $F^m(x) \sim \exp(-M_n)$. Then, the following Corollary 1 holds, which is a direct consequence of Theorem 2.2 in \cite{dombry2019maximum} providing asymptotic properties of the maximum likelihood estimator $\widehat{\bm{\gamma}}$. The symbol $\sim$ means $p_n \sim q_n ~\Longleftrightarrow~ p_n / q_n = 1+ o(1)$. Hereafter, all asymptotic notations refer to $n \to \infty$.

\begin{corol}
Under some conditions $(F^m(x) - G_{\widehat{\bm{\gamma}}}(x))$ converges with the rate the larger of $F^m(x) - G_{{\bm{\gamma}}}(x)$ and
$$\sqrt{\frac{m}{n}} M_n \exp(-M_n) \times
\begin{dcases}
\lambda_n (M_n^{\gamma} + \ln M_n) + \ln M_n ~~ &\text{for} ~~ M_n \to 0\\
\lambda_n + 1~~ &\text{for} ~~ M_n =O(1) \\
\lambda_n (M_n^{\gamma} + \ln M_n) + M_n^{\gamma} ~~ &\text{for} ~~ M_n \to \infty,
\end{dcases}
$$
where
$$\lambda_n := \sqrt{m} \times
\begin{dcases}
 m^{-\beta} ~~~ &\text{for the Hall class} \\
(\ln m)^{-1} ~~~ &\text{for the Weibull class} \\
m^{\sigma} ~~~ &\text{for the bounded class}.
\end{dcases}
$$
\end{corol} 
The details including the exact condition and the proof are given in the Appendices.

The approximation $F^m(x) - G_{{\bm{\gamma}}}(x)$ can be quite poor, where the normal distribution cases are well known and studied (see \citealt{hall1979rate}; \citealt{hall1980estimating}). To improve the approximation error, the Edgeworth expansion shown in Theorem 3.1 in \cite{haan1996second} may be helpful. Increasing the model degree is another way to reduce the approximation error (e.g. applying the flexible $q$-generalized extreme value distribution and exponential generalization proposed by \citealt{nascimento2015extended} and \citealt{serge2018generalized}, respectively). While these approaches can be applied to error reduction, the approximation cannot be improved significantly.

It is obvious that the fitting to SMD is unsuitable for distributions not belonging to any maximum domain of attractions. Besides, $F^m$ can be far from $G_{{\bm{\gamma}}}$ for not large $m$ unless $F$ itself is a generalized extreme value distribution. Motivated by the above facts, we change our way from extreme value based one to nonparametric one. Nonparametric approaches have been applied to extreme value analysis in different settings especially for describing dependence structure of variables (\citealt{beirlant2003regression}, \citealt{beirlant2004local}, \citealt{wang2009tail}, \citealt{beirlant2004nonparametric},  \citealt{yoshida2021nonparametric}).

We consider the plug-in type of estimator defined as a nonparametric distribution estimator to the power of $m$. We call the nonparametric estimator NE. The `parametrically' fitting estimator $G_{\widehat{\bm{\gamma}}}$ is called PE. NE does not need model specification and can be consistent even for cases that the Fisher-Tippett-Gnedenko theorem breaks as long as a condition holds including the twice-differentiability (the details introduced in Section 2). Moreover, NE is expected to be better than PE at least for relatively small $m$ due to the wide applicability. It follows that NE seems to be a good candidate for SMD. The three types of distribution are assumed in order to obtain the explicit forms of asymptotic properties in this study; however, we note that NE does not need such a model specification.

Asymptotic properties of NE are provided in Section 2. The results of theoretical and numerical studies comparing the accuracy of the two different approaches are given in Section 3. It is demonstrated that numerical performances heavily depend on not only $m$ but also the extreme value index $\gamma$, where they are comparable as a whole. Two real case studies and a discussion are presented in Sections 4 and 5, respectively. The details of PE and the proofs are given in the Appendices. 

\section{Kernel-type estimation of SMD}

This study employs the kernel distribution estimator as the nonparametric distribution estimator. The kernel distribution estimator is given by
$$\widehat{F}(x)=\frac{1}{n} \sum_{i=1}^n W\left(\frac{x -X_i}{h} \right),$$
where $W$ is a smooth cumulative distribution function. Let $w$ be the symmetric and bounded density function, whose support is bounded for distributions with bounded support. $h$ is the bandwidth that satisfies $h \rightarrow 0$ and for the bounded class of distributions $h (x^* -x)^{-1} \to 0$. 

The following theorem on the mean squared error (MSE) and the asymptotic normality of NE holds for the three types of distributions.

\begin{theorem}
Suppose $F$ is continuously twice differentiable at $x$ and $\int z^2 w(z) {\rm d}z <\infty$. Then,
\begin{align*}
\exp(2M_n) \mathbb{E}[(F^m(x) -\widehat{F}^m(x))^2] \sim \Biggm\{ M_n h^{2} \frac{\xi_n}{2} \int z^2 w(z) {\rm d}z  \Biggm\}^2 + \frac{m}{n} M_n,
\end{align*}
where
$$
\xi_n:=\begin{dcases}
\alpha(\alpha+1) x^{-2}~~~ &\text{for the Hall class}\\
\kappa^2 C^2 x^{2\kappa-2} ~~~ &\text{for the Weibull class}\\
\mu(\mu+1) (x^* -x)^{-2} ~~~ &\text{for the bounded class}.
\end{dcases}
$$
Furthermore, 
$$\exp(M_n) \left( \frac{m}{n} M_n \right)^{-1/2} \{F^m(x) -\widehat{F}^m(x) + \exp(-M_n) M_n h^{2} \frac{\xi_n}{2} \int z^2 w(z) {\rm d}z\}$$
converges in distribution to the standard normal distribution if the MSE converges to zero.
\end{theorem}

\begin{remark}
The condition that $F$ is continuously twice differentiable is broken for the bounded class with $\mu\le -2$ or with $\mu=-1$ and $\sigma \le -1$.
\end{remark}

\begin{remark}
For the bounded class of distributions $x$ satisfying $\lim_{n\to \infty} h(x^* -x)^{-1} >0$ is called a boundary point in the kernel distribution estimation and the convergence rate of $\widehat{F}^m(x)$ changes (see the proof of Theorem 1).
\end{remark}

The following Corollary 2 on the optimal convergence rate follows from Theorem 1.
\begin{corol}
Suppose $|\int z W(z) w(z) {\rm d}z| < \infty$ and $\xi_n^{-2} \psi_n n^{-1} \to 0$ hold, where 
$$
\psi_n:=\begin{dcases}
A^{-1} \alpha x^{\alpha-1} ~~~ &\text{for the Hall class}\\
\kappa C x^{\kappa-1} \exp(C x^{\kappa})~~~ &\text{for the Weibull class}\\
-D^{-1} \mu (x^* -x)^{\mu -1} ~~~ &\text{for the bounded class}.
\end{dcases}
$$ 
Under the assumptions of Theorem 1, the asymptotically optimal bandwidth is
$$h= \left(2 \xi_n^{-2} \psi_n n^{-1} \frac{\int z W(z) w(z) {\rm d}z}{\left(\int z^2 w(z) {\rm d}z\right)^2}\right)^{1/3}.$$ 
If the optimal bandwidth converges to 0, $F^m(x) -\widehat{F}^m(x)$ with the optimal bandwidth is asymptotically nondegenerate normal with the asymptotic bias 
$$\nu_0 \exp(-M_n) M_n n^{-2/3} \xi_n^{-1/3} \psi_n^{2/3},$$
where $\nu_0:= \left(2 \int z^2 w(z) {\rm d}z\right)^{-1/3} \left(\int z W(z) w(z) {\rm d}z\right)^{2/3}$.
\end{corol}

\begin{remark}
The convergence condition of the optimal bandwidth is restrictive. However, any $h$ converging to 0 ensures the consistency and can be the candidate in practice as long as Theorem 1, $M_n \xi_n \to 0$ and $|\int z W(z) w(z) {\rm d}z| < \infty$ hold.
\end{remark}

\noindent
The following Corollary 3 states the special case $M_n=O(1)$ of Corollary 2.
\begin{corol}
Suppose Corollary 2 holds and $\exists\delta>0$ s.t. $M_n \to \delta$. Then, the asymptotically optimal bandwidth is
\begin{align*}
h=& \left(2 \delta^{-1-3\gamma} m^{1+3\gamma} n^{-1} \frac{\int z W(z) w(z) {\rm d}z}{\left(\int z^2 w(z) {\rm d}z\right)^2}\right)^{1/3} \\
& \times
\begin{dcases}
A^{\gamma} \{\sqrt{\alpha}(\alpha+1)\}^{-2/3} ~~~ &\text{for the Hall class}\\
(\kappa C)^{-1} (C^{-1} \ln (m/\delta))^{-\theta} ~~~ &\text{for the Weibull class}\\
 D^{\gamma} \{\sqrt{-\mu}(-\mu+1)\}^{-2/3} ~~~ &\text{for the bounded class}.
\end{dcases}
\end{align*}
If the optimal bandwidth converges to 0, $F^m(x) -\widehat{F}^m(x)$ with the optimal bandwidth has the asymptotic bias
\begin{align*}
\nu_0 \exp(-\delta) \delta^{1/3} m^{2/3} n^{-2/3} \times
\begin{dcases}
(\alpha (\alpha+1)^{-1})^{1/3} ~~~ &\text{for the Hall class}\\
1 ~~~ &\text{for the Weibull class}\\
(\mu (\mu+1)^{-1})^{1/3} ~~~ &\text{for the bounded class}.
\end{dcases}
\end{align*}
\end{corol}
The requirement for the bounded class of distributions $h (x^* - x) \to 0$ means $m=o(n)$ (see Remark 3).

Note that the asymptotically optimal bandwidth possibly does not converge to 0. In such cases the asymptotic bias does not converge; however, any bandwidth converging to 0 makes NE consistent. Hence, the optimal bandwidth subject to $h \to 0$ is considered to be slower than any polynomial (e.g. $h=(\ln n)^{-1}$). The asymptotic bias converges to 0 with the rate $M_n h^2 \xi_n$, and it is of order $m^{-2\gamma} h^2$.

The proof of Theorem 1 shows that the MSE asymptotically coincides with $m^2$ times that of the kernel distribution estimator $\widehat{F}(x)$, which means that the asymptotically optimal bandwidth of NE is independent of $m$, as shown in Corollary 2. It follows that minimizing the MSE of NE is asymptotically equivalent to minimizing that of $\widehat{F}(x)$.

The optimal bandwidth for the pointwise estimation needs to estimate derivatives of the distribution $F$ at $x$, which means that the extreme value index needs to be estimated. Therefore, the optimal bandwidth estimator (i.e., NE itself) must include an extreme value index estimator, which is considered to affects the convergence rate of NE substantially. However, unlike the MLE $\widehat{\bm{\gamma}}$ of PE, there are no restrictions on the block size (see the Appendices) in the  bandwidth selection. It is possible to choose the optimal the block size for the extreme value index estimation, and NE with the bandwidth estimator possibly theoretically outperforms PE as a result.

Section 3 reports the results of numerical studies on the mean integrated squared error (MISE) of NE with the bandwidth estimator of the kernel distribution estimator $\widehat{F}$ proposed by \cite{altman1995bandwidth}. The \texttt{kerdiest} package in the R software environment provides the \texttt{aLbw} function for calculating the plug-in value. Although the bandwidth estimator asymptotically minimizes not the MSE of NE but the MISE of $\widehat{F}$, we report it as the case study for the purpose of measuring the applicability. It is considered that the kernel function selection does not significantly affect the numerical accuracy (\cite{lejeune1992smooth}).

\section{Comparative study on the SMD estimators}

Suppose $\exists\delta>0$ s.t. $M_n \equiv \delta$ throughout in this section. The MSE of PE $\mathbb{E}[(F^m(x) - G_{\widehat{\bm{\gamma}}}(x))^2]$ converges with the rate $\{(m/n) + m^{-2\beta\gamma} + m^{-2}\}$ for the Hall class and $\{(m/n) +  m^{-2\sigma\gamma} + m^{-2}\}$ for the bounded class. For the Weibull class PE is found to be inconsistent by applying $M_n \equiv \delta$ to Corollary 1. The MSE of NE with the optimal bandwidth is of order $(m/n)$ if the optimal bandwidth converges to zero i.e. $m^{1+3\gamma} n^{-1} \to 0$. The underlying distributions $F$ are supposed to be Pareto distributions, T distributions, Burr distributions, Fr\'{e}chet distributions, Weibull distributions or inverse Burr distributions. 

By simulating the following MISE of PE $G_{\widehat{\bm{\gamma}}}$
$$L_m^{-1} \int_{Q_ m(0.1)}^{Q_ m(0.9)} \left(G_{\widehat{\bm{\gamma}}}(x) -F^m(x)\right)^2 {\rm d}x$$ 
and that of NE $\widehat{F}^m$, we studied the accuracy in finite-sample cases, where $L_m := Q_ m(0.9) - Q_ m(0.1)$ and $Q_ m(q)$ denotes the $q$th quantile of SMD. We suppose $Q_ m(q) ~ (0.1 \le q \le 0.9)$ means $M_n=O(1)$. The parameters of the underlying distributions, the  convergence rates of the MSE without terms slower than any polynomial, and length $L_m$ for $n=2^{12}$ are summarized in Table 1. We set $m =n^{1/4}$, $m=n^{1/2}$, and $m=n^{3/4}$. Table 1 shows that the slowness of the convergence rate of PE for large $\alpha$. A hyphen in Table 1 means the underlying distribution breaks the assumption of this study. We see the assumptions of both of PE and NE are restrictive especially for the bounded class, but we report all the numerical results for the sake of future research. 

We simulated the MISE values 10000 times, where Tables 2--5 show the mean values and their standard deviations (sd). The sample size $n$ was $2^{8}$ or $2^{12}$. The bandwidth was estimated by \cite{altman1995bandwidth}'s method. The kernel function was the Epanechnikov for the inverse Burr distributions and the Gaussian for the other distributions. The block size in the maximum likelihood method (see the Appendices) was chosen as $(\ln m)^2$ for the Weibull distributions so that $\lambda_n \to 1$ and as same as $m$ for the other distributions. 

Comparing PE and NE for the Pareto distributions with the shape parameters $\ell=1/2, 1, 3, 10$, where the corresponding parameters $\alpha=\ell$ and $\beta=1$, we see that NE surpassed PE when at least one of $m$ and $\ell$ is large. In the other cases, the MISE values of PE were smaller. For the T distributions, which are only different from the Pareto distributions in the sense of the second-order parameter $\beta$, both PE and NE have similar numerical properties. The accuracy of NE for the T distribution with $\ell=1/2$ is worse than that for the Pareto distribution with $\ell=1/2$. The Burr distributions have two parameters, $c$ and $\ell$, where $\alpha=c\ell$ and $\beta=c$. Table 3 shows that for the Burr distributions NE tends to outperform PE when at least one of $m$ and $c \times \ell$ is large. The Fr\'{e}chet distributions have three parameters: the location, scale, and shape parameters. In this study, the location parameter was set to 0, the scale parameter was set to 1, and the shape parameter was set to $5, 2, 1, 1/2, 1/4$. In the Fr\'{e}chet cases, NE was better for small $\gamma$. Contrary to the theoretical results shown in Table 1 the MISE of NE sometimes decreases with $m$ for the heavy-tailed Pereto, T, Burr, Fr\'{e}chet cases. This may come from the numerical property that the kernel method works poorly for heavy-tailed distributions (in density estimation see \citealt{maiboroda2004estimation}).

In the Weibull cases with $D=1$ and $\kappa=1/2, 1, 3, 10$, which means $\gamma=0$ and the second-order parameter $\rho=0$, NE surpassed PE in every case. PE is not consistent for the inverse Burr distribution with $\mu \ge -2$. The distribution with the non-integer $\mu$ more than $-2$ changes the usual convergence rate of NE (see Remark 2). Though NE also requires $h (x^* -x)^{-1} \to 0$, it holds for $M_n \equiv \delta$ since $h (x^* -x)^{-1} =O((m/n)^{1/3})$. The MISE values of PE are large for $\mu$ close to zero, and those of NE are large especially for $c=1$. 

The numerical study clarified the large degree of misspecification of PE for $\gamma$ close to zero. In this case NE numerically outperformed PE not only for small $m$ but also large $m$. For $\gamma$ far from zero NE was not accurate in general; however, the MISE values of NE in such cases sometimes get smaller as $m$ becomes large. Therefore, the study demonstrated the possibility of numerical superiority of NE for large $m$ and that the relative numerical performance in finite-sample cases heavily depend on not only $m$ but also $\gamma$. Summarizing the above, the accuracy of the two approaches heavily depends on cases, but NE is comparable to PE on the whole. 

\begin{table}
\caption{The polynomial convergence rates of the MSE of the estimators and the lengths}{\fontsize{8pt}{8pt}\selectfont
$$\begin{tabu}[c]{|ccc||ccc|ccc|ccc|}
  \hline
 \multicolumn{3}{|c||}{{\rm Pareto}} & \multicolumn{3}{c|}{m=n^{1/4}} & \multicolumn{3}{c|}{m=n^{1/2}} & \multicolumn{3}{c|}{m=n^{3/4}} \\
\ell & \alpha & \beta & {\rm PE} & {\rm NE} & (L_m) & {\rm PE} & {\rm NE} & (L_m) & {\rm PE} & {\rm NE} & (L_m) \\ \hline

 1/2 & 1/2 & 1 & -1/2 & \text{--} & (6 \times 10^3) & -1/2 & \text{--} & (4 \times 10^5) & -1/4 & \text{--} & (2 \times 10^{7})  \\ 
 1 & 1 & 1 & -1/2 & \text{--} & (70) & -1/2 & \text{--} & (600) & -1/4 & \text{--} & (5 \times 10^3) \\ 
 3 & 3 & 1 & -1/6 & -3/4 & (3) & -1/3 & \text{--} & (5) & -1/4 & \text{--} & (10) \\ 
 10 & 10 & 1 & -1/20 & -3/4 & (0.4) & -1/10 & -1/2 & (0.5) & -3/20 & -1/4 & (0.6)  \\ 
 
\hline\hline
 \multicolumn{3}{|c||}{{\rm T}} & \multicolumn{3}{c|}{m=n^{1/4}} & \multicolumn{3}{c|}{m=n^{1/2}} & \multicolumn{3}{c|}{m=n^{3/4}} \\
\ell & \alpha & \beta & {\rm PE} & {\rm NE} & (L_m) & {\rm PE} & {\rm NE} & (L_m) & {\rm PE} & {\rm NE} & (L_m) \\ \hline

 1/2 & 1/2 & 2 & -1/2 & \text{--} & (600) & -1/2 & \text{--} & (4 \times 10^4) & -1/4 & \text{--} & (2 \times 10^{6}) \\ 
 1 & 1 & 2 & -1/2 & \text{--} & (20) & -1/2 & \text{--} & (200) & -1/4 & \text{--} & (10^3) \\ 
 3 & 3 & 2 & -1/3 & -3/4 & (3) & -1/2 & \text{--} & (6) & -1/4 & \text{--} & (10) \\ 
 10 & 10 & 2 & -1/10 & -3/4 & (2) & -1/5 & -1/2 & (2) & -1/4 & -1/4 & (2) \\ 

 \hline\hline
  \multicolumn{3}{|c||}{{\rm Burr}} & \multicolumn{3}{c|}{m=n^{1/4}} & \multicolumn{3}{c|}{m=n^{1/2}} & \multicolumn{3}{c|}{m=n^{3/4}} \\
c,\ell & \alpha & \beta & {\rm PE} & {\rm NE} & (L_m) & {\rm PE} & {\rm NE} & (L_m) & {\rm PE} & {\rm NE} & (L_m) \\ \hline

 1/2,1/2 & 1/4 & 1/2 & -1/2 & \text{--} & (3 \times 10^7) & -1/2 & \text{--} & (10^{11}) & -1/4 & \text{--} & (6 \times 10^{14})  \\ 
 1,1/2 & 1/2 & 1 & -1/2 & \text{--} & (6 \times 10^3) & -1/2 & \text{--} & (4 \times 10^{5}) & -1/4 & \text{--} & (2 \times 10^{7})  \\ 
 3,1/2 & 3/2 & 3 & -1/2 & -3/4 & (10) & -1/2 & \text{--} & (50) & -1/4 & \text{--} & (200)  \\ 
 1/2,1 & 1/2 & 1/2 & -1/2 & \text{--} & (6 \times 10^3) & -1/2 & \text{--} & (4 \times 10^5) & -1/4 & \text{--} & (2 \times 10^{7}) \\ 
 1,1 & 1 & 1 & -1/2 & \text{--} & (70) & -1/2 & \text{--} & (600) & -1/4 & \text{--} & (5 \times 10^3) \\  
 3,1 & 3 & 3 & -1/2 & -3/4 & (3) & -1/2 & \text{--} & (5) & -1/4 & \text{--} & (10) \\ 
 1/2,3 & 3/2 & 1/2 & -1/6 & -3/4 & (20) & -1/3 & \text{--} & (60) & -1/4 & \text{--} & (300) \\ 
 1,3 & 3 & 1 & -1/6 & -3/4 & (3) & -1/3 & \text{--} & (5) & -1/4 & \text{--} & (10)  \\ 
 3,3 & 9 & 3 & -1/6 & -3/4 & (0.6) & -1/3 & -1/2 & (0.7) & -1/4 & \text{--} & (0.8)  \\ 
 
  \hline\hline
 \multicolumn{3}{|c||}{{\rm Fr{e}chet}} & \multicolumn{3}{c|}{m=n^{1/4}} & \multicolumn{3}{c|}{m=n^{1/2}} & \multicolumn{3}{c|}{m=n^{3/4}} \\
\gamma & \alpha & \beta & {\rm PE} & {\rm NE} & (L_m) & {\rm PE} & {\rm NE} & (L_m) & {\rm PE} & {\rm NE} & (L_m) \\ \hline

 5 & 1/5 & 1/5 & \text{--} & \text{--} & (5 \times 10^{8}) & \text{--} & \text{--} & (2 \times 10^{13}) & \text{--} & \text{--} & (5 \times 10^{17})\\ 
 2 & 1/2 & 1/2 & -1/2 & \text{--} & (3 \times 10^{3}) & -1/2 & \text{--} & (2 \times 10^{5}) & -1/4 & \text{--} & (10^{7})\\  
 1 & 1 & 1 & -1/2 & \text{--} & (70) & -1/2 & \text{--} & (600) & -1/4 & \text{--} & (5 \times 10^3) \\ 
 1/2 & 2 & 1 & -1/4 & -3/4 & (10) & -1/2 & \text{--} & (40) & -1/4 & \text{--} & (100) \\ 
 1/4 & 4 & 1 & -1/8 & -3/4 & (6) & -1/4 & -1/2 & (10) & -1/4 & \text{--} & (20) \\
 
 \hline\hline
 \multicolumn{3}{|c||}{{\rm Weibull}} & \multicolumn{3}{c|}{m=n^{1/4}} & \multicolumn{3}{c|}{m=n^{1/2}} & \multicolumn{3}{c|}{m=n^{3/4}} \\
\kappa& \gamma & \rho & {\rm PE} & {\rm NE} & (L_m) & {\rm PE} & {\rm NE} & (L_m) & {\rm PE} & {\rm NE} & (L_m) \\ \hline

1/2 & 0 & 0 & \text{--} & -3/4 & (20) & \text{--} & -1/2 & (30) & \text{--} & -1/4 & (40) \\
1 & 0 & 0 & \text{--} & -3/4 & (3) & \text{--} & -1/2 & (3) & \text{--} & -1/4 & (3) \\
3 & 0 & 0 & \text{--} & -3/4 & (0.5) & \text{--} & -1/2 & (0.4) & \text{--} & -1/4 & (0.3) \\
10 & 0 & 0 & \text{--} & -3/4 & (0.1) & \text{--} & -1/2 & (0.08) & \text{--} & -1/4 & (0.05) \\

 \hline\hline
  \multicolumn{3}{|c||}{{\rm inv. Burr}} & \multicolumn{3}{c|}{m=n^{1/4}} & \multicolumn{3}{c|}{m=n^{1/2}} & \multicolumn{3}{c|}{m=n^{3/4}} \\
c,\ell & \mu & \sigma & {\rm PE} & {\rm NE} & (L_m) & {\rm PE} & {\rm NE} & (L_m) & {\rm PE} & {\rm NE} & (L_m) \\ \hline

3, 2 & -6 & -2 & -1/6 & -3/4 & (0.7) & -1/3 & -1/2 & (0.3)  & -1/4 & -1/4 & (0.2) \\
1, 2 & -2 & -2 & \text{--} & -3/4 & (0.5) & \text{--} & -1/2 & (0.2)  & \text{--} & -1/4 & (0.05) \\
1/2, 2 & -1 & -2 & \text{--} & -3/4 & (0.2) & \text{--} & -1/2 & (0.03)  & \text{--} & -1/4 & (4 \times 10^{-3}) \\
 3, 1 & -3 & -1 & -1/6 & -3/4 & (1) & -1/3 & -1/2 & (0.4)  & -1/4 & -1/4 & (0.1) \\
1, 1 & -1 & -1 & \text{--} & -3/4 & (0.3) & \text{--} & -1/2 & (0.03)  & \text{--} & -1/4 & (4 \times 10^{-3}) \\
1/2,1 & -1/2 & -1 & \text{--} & \text{--} & (0.07) & \text{--} & \text{--} & (10^{-3})  & \text{--} & \text{--} & (2 \times 10^{-5}) \\
 3, 1/3 & -1 & -1/3 & \text{--} & \text{--} & (5) & \text{--} & \text{--} & (0.1)  & \text{--} & \text{--} & (7 \times 10^{-3}) \\
 1,1/3 & -1/3 & -1/3 & \text{--} & \text{--} & (0.04) & \text{--} & \text{--} & (5 \times 10^{-5})  & \text{--} & \text{--} & (9 \times 10^{-8}) \\
1/2, 1/3 & -1/6 & -1/3 & \text{--} & \text{--} & (3 \times 10^{-4}) & \text{--} & \text{--} & (2 \times 10^{-9}) & \text{--} & \text{--} & (8 \times 10^{-15})\\

 \hline
\end{tabu}$$
}
\end{table}

\begin{table}
\caption{The scaled MISE values ($\times 100$) and sd values ($\times 100$) for the estimators}
{\fontsize{10pt}{9pt}\selectfont
$$\begin{tabu}[c]{|c||cccc|cccc|}
   \hline
    {\rm Pareto} & \multicolumn{4}{c|}{n=2^{8}} & \multicolumn{4}{c|}{n=2^{12}} \\ 
     & {\rm PE} & {\rm sd} & {\rm NE} & {\rm sd} & {\rm PE} & {\rm sd} & {\rm NE} & {\rm sd} \\  \hline\hline 
    
\ell & \multicolumn{4}{c|}{m=2^{2}} & \multicolumn{4}{c|}{m=2^{3}} \\ \hline
1/2 & 0.328 & 1.950 & 34.356 & 23.599 & 0.218 & 1.291 & 60.708 & 16.017 \\
  1 & 0.168 & 0.204 & 4.936 & 11.175 & 0.024 & 0.263 & 4.849 & 14.329 \\
  3 & 0.201 & 0.233 & 0.208 & 0.213 & 0.025 & 0.030 & 0.021 & 0.019 \\
  10 & 0.214 & 0.250 & 0.208 & 0.213 & 0.026 & 0.031 & 0.026 & 0.027 \\
   
   \hline\hline
\ell & \multicolumn{4}{c|}{m=2^{4}} & \multicolumn{4}{c|}{m=2^{6}} \\ \hline
1/2 & 4.524 & 7.009 & 22.191 & 28.085 & 0.473 & 2.091 & 23.318 & 27.790 \\
  1 & 1.046 & 2.385 & 3.820 & 11.028 & 0.241 & 0.903 & 2.516 & 10.588 \\
  3 & 0.870 & 0.984 & 0.808 & 0.820 & 0.200 & 0.241 & 0.188 & 0.194 \\
  10 & 0.918 & 1.034 & 0.788 & 0.805 & 0.203 & 0.241 & 0.177 & 0.142 \\
   
      \hline\hline
\ell & \multicolumn{4}{c|}{m=2^{6}} & \multicolumn{4}{c|}{m=2^{9}} \\ \hline
1/2 & 12.469 & 12.822 & 11.023 & 20.566 & 9.442 & 10.118 & 5.679 & 15.210 \\
  1  & 8.539 & 6.504 & 4.119 & 7.716 & 5.829 & 5.768 & 2.883 & 8.423 \\
  3 & 5.564 & 4.703 & 3.086 & 2.708 & 2.061 & 2.202 & 1.634 & 1.613 \\
  10 & 5.191 & 4.667 & 3.047 & 2.867 & 2.186 & 2.311 & 1.659 & 1.515 \\
   
    \hline\hline
{\rm T} & \multicolumn{4}{c|}{n=2^{8}} & \multicolumn{4}{c|}{n=2^{12}} \\
 & {\rm PE} & {\rm sd} & {\rm NE} & {\rm sd} & {\rm PE} & {\rm sd} & {\rm NE} & {\rm sd} \\  \hline\hline 
 \ell & \multicolumn{4}{c|}{m=2^{2}} & \multicolumn{4}{c|}{m=2^{3}} \\ \hline

1/2 & 0.265 & 1.083 & 44.619 & 19.123 & 0.112 & 0.990 & 67.907 & 5.001 \\
  1 & 0.237 & 0.276 & 5.920 & 11.863 & 0.030 & 0.031 & 6.885 & 15.454 \\
  3 & 0.223 & 0.246 & 0.191 & 0.214 & 0.033 & 0.035 & 0.028 & 0.031 \\
  10 & 0.208 & 0.240 & 0.184 & 0.214 & 0.029 & 0.030 & 0.026 & 0.037 \\

\hline\hline
 \ell & \multicolumn{4}{c|}{m=2^{4}} & \multicolumn{4}{c|}{m=2^{6}} \\ \hline
1/2  & 4.064 & 6.836 & 33.468 & 29.871 & 0.385 & 1.903 & 45.261 & 28.737 \\
  1 & 0.875 & 1.854 & 4.981 & 12.725 & 0.171 & 0.231 & 2.295 & 8.902 \\
  3  & 0.896 & 0.985 & 0.720 & 0.759 & 0.200 & 0.237 & 0.235 & 0.245 \\
  10 & 0.934 & 1.047 & 0.644 & 0.730 & 0.207 & 0.245 & 0.195 & 0.174 \\

\hline\hline
\ell & \multicolumn{4}{c|}{m=2^{6}} & \multicolumn{4}{c|}{m=2^{9}} \\ \hline
1/2  & 12.359 & 12.221 & 16.353 & 24.658 & 9.819 & 10.261 & 13.869 & 23.565 \\
  1  & 8.346 & 6.477 & 4.633 & 8.997 & 5.371 & 5.678 & 1.594 & 1.438 \\
  3  & 5.517 & 4.748 & 2.846 & 2.631 & 2.097 & 2.172 & 1.708 & 1.296 \\
  10  & 5.335 & 4.745 & 2.315 & 2.439 & 2.226 & 2.305 & 1.422 & 1.589 \\
 \hline
\end{tabu}$$
}
\end{table}

\begin{table}
\caption{The scaled MISE values ($\times 100$) and sd values ($\times 100$) for the estimators}
{\fontsize{10pt}{9pt}\selectfont
$$\begin{tabu}[c]{|c||cccc|cccc|}
   \hline
    {\rm Burr} & \multicolumn{4}{c|}{n=2^{8}} & \multicolumn{4}{c|}{n=2^{12}} \\
  & {\rm PE} & {\rm sd} & {\rm NE} & {\rm sd} & {\rm PE} & {\rm sd} & {\rm NE} & {\rm sd} \\  \hline\hline 
    
 c,\ell & \multicolumn{4}{c|}{m=2^{2}} & \multicolumn{4}{c|}{m=2^{3}}  \\ \hline
1/2,1/2 & 0.801 & 3.808 & 60.715 & 14.740 & 5.671 & 10.605 & 75.900 & 0.000 \\
  1,1/2 & 0.351 & 2.111 & 34.113 & 23.627 & 0.207 & 1.368 & 59.039 & 18.150 \\
    3,1/2 & 0.184 & 0.224 & 0.621 & 2.826 & 0.022 & 0.027 & 0.036 & 0.116 \\
  1/2,1 & 0.366 & 1.978 & 35.100 & 23.664 & 0.263 & 1.424 & 59.045 & 19.638 \\
    1,1 & 0.167 & 0.202& 5.110 & 11.410 & 0.021 & 0.095 & 3.148 & 10.356 \\
    3,1 & 0.193 & 0.225 & 0.196 & 0.210 & 0.024 & 0.029 & 0.024 & 0.024 \\
    1/2,3 & 0.173 & 0.201 & 1.842 & 5.574 & 0.028 & 0.030 & 0.403 & 3.264 \\
     1,3 & 0.201 & 0.237 & 0.208 & 0.219 & 0.024 & 0.029 & 0.033 & 0.038 \\
       3,3 & 0.206 & 0.242 & 0.193 & 0.225 & 0.027 & 0.030 & 0.024 & 0.028 \\
       
 \hline\hline
 c,\ell & \multicolumn{4}{c|}{m=2^{4}} & \multicolumn{4}{c|}{m=2^{6}}  \\ \hline
1/2,1/2 & 3.860 & 8.797 & 45.726 & 33.915 & 1.536 & 4.437 & 65.212 & 23.362 \\
  1,1/2 & 4.691 & 7.131 & 21.580 & 27.710 & 0.367 & 2.068 & 24.791 & 29.590 \\
  3,1/2 & 0.794 & 0.984 & 1.050 & 3.091 & 0.177 & 0.207 & 0.177 & 0.190 \\
  1/2,1 & 4.622 & 7.022 & 22.202 & 27.983 & 0.508 & 2.332 & 23.275 & 28.204 \\
  1,1 & 1.082 & 2.541 & 3.945 & 11.081 & 0.237 & 0.897 & 1.199 & 6.698 \\
  3,1 & 0.877 & 0.992 & 0.782 & 0.838 & 0.197 & 0.228 & 0.210 & 0.216 \\
  1/2,3  & 0.798 & 1.179 & 1.552 & 5.095 & 0.181 & 0.216 & 0.284 & 0.418 \\
  1,3 & 0.870 & 0.969 & 0.786 & 0.811 & 0.198 & 0.235 & 0.198 & 0.250 \\
  3,3  & 0.963 & 1.121 & 0.650 & 0.705 & 0.203 & 0.241 & 0.183 & 0.169 \\
  
 \hline\hline
 c,\ell & \multicolumn{4}{c|}{m=2^{6}} & \multicolumn{4}{c|}{m=2^{9}}  \\ \hline
1/2,1/2 & 16.478 & 21.085 & 22.598 & 31.978 & 12.798 & 18.836 &  22.562 & 32.222 \\ 
  1,1/2  & 12.421 & 12.748 & 10.090 & 19.481 & 10.290 & 10.129 & 5.744 & 15.562 \\ 
  3,1/2  & 6.928 & 5.286 & 3.081 & 3.367 & 3.055 & 3.750 & 1.422 & 1.397 \\ 
  1/2,1  & 12.450 & 12.775 & 10.397 & 19.295 & 9.585 & 10.358 & 6.159 & 16.112 \\ 
  1,1  & 8.510 & 6.564 & 4.115 & 7.731 & 5.758 & 5.728 & 1.572 & 1.226 \\ 
   3,1  & 5.550 & 4.713 & 3.104 & 2.725 & 2.084 & 2.195 & 1.414 & 1.155 \\ 
   1/2,3  & 7.120 & 5.446 & 3.330 & 4.156 & 3.025 & 3.721 & 1.513 & 1.415 \\ 
   1,3  & 5.616 & 4.840 & 3.051 & 2.523 & 2.099 & 2.250 & 1.768 & 1.659 \\ 
   3,3  & 5.416 & 4.812 & 2.496 & 2.480 & 2.209 & 2.336 & 1.572 & 1.237 \\ 
   
 \hline
\end{tabu}$$
}
\end{table}

\begin{table}
\caption{The scaled MISE values ($\times 100$) and sd values ($\times 100$) for the estimators}
{\fontsize{10pt}{9pt}\selectfont
$$\begin{tabu}[c]{|c||cccc|cccc|}
   \hline
  {\rm Frechet} & \multicolumn{4}{c|}{n=2^{8}} & \multicolumn{4}{c|}{n=2^{12}} \\
  & {\rm PE} & {\rm sd} & {\rm NE} & {\rm sd} & {\rm PE} & {\rm sd} & {\rm NE} & {\rm sd} \\  \hline\hline 
    
\gamma & \multicolumn{4}{c|}{m=2^{2}} & \multicolumn{4}{c|}{m=2^{3}}  \\ \hline
5 & 0.203 & 0.239 & 63.838 & 12.297 & 0.024 & 0.029 & 77.1 & 0 \\
2 & 0.191 & 0.225 & 33.904 & 23.627 & 0.023 & 0.028 & 56.019 & 20.548 \\
    1 & 0.166 & 0.204 & 5.168 & 11.675 & 0.023 & 0.196 & 3.839 & 10.325 \\
  1/2 & 0.282 & 1.661 & 0.307 & 1.581 & 0.225 & 1.408 & 0.028 & 0.029 \\
1/4 & 1.288 & 6.490 & 0.205 & 0.216& 11.206 & 16.847 & 0.026 & 0.026 \\
        
 \hline\hline
\gamma & \multicolumn{4}{c|}{m=2^{4}} & \multicolumn{4}{c|}{m=2^{6}}  \\ \hline
5  & 0.896 & 1.027 & 51.211 & 34.094 & 0.197 & 0.234 & 69.953 & 19.804 \\
  2  & 0.845 & 0.956 & 21.994 & 28.278 & 0.194 & 0.226 & 24.115 & 30.162 \\
  1  & 1.102 &  2.550 & 4.191 & 11.844 & 0.242 & 0.918 & 0.816 & 4.627 \\
  1/2  & 4.428 &  6.920 & 0.954 & 2.080 & 0.459 & 2.083 & 0.221 & 0.237 \\
  1/4 & 3.965 & 10.296 & 0.772 & 0.761 & 2.602 & 8.544 & 0.202 & 0.205 \\
 \hline\hline
\gamma & \multicolumn{4}{c|}{m=2^{6}} & \multicolumn{4}{c|}{m=2^{9}} \\ \hline
5  & 5.849 & 4.919 & 26.030 & 34.162 & 2.090 & 2.172 & 28.942 & 34.469 \\
2  & 6.356 & 4.762 & 10.641 & 19.907 & 2.394 & 2.831 & 9.425 & 20.763 \\
  1  & 8.568 & 6.565 & 4.275 & 7.958 & 5.823 & 5.746 & 1.527 & 1.265 \\
  1/2  & 12.802 & 12.850 & 3.323 & 3.182 & 9.791 & 10.531 & 1.567 & 1.389 \\
  1/4  & 18.141 & 23.672 & 3.048 & 2.662 & 14.080 & 21.917 & 1.614 & 1.510 \\
 \hline
   \hline
    {\rm Weibull} & \multicolumn{4}{c|}{n=2^{8}} & \multicolumn{4}{c|}{n=2^{12}} \\
 & {\rm PE} & {\rm sd} & {\rm NE} & {\rm sd} & {\rm PE} & {\rm sd} & {\rm NE} & {\rm sd} \\  \hline\hline 
    
\kappa & \multicolumn{4}{c|}{m=2^{2}} & \multicolumn{4}{c|}{m=2^{3}} \\ \hline
1/2 & 2.768 & 0.683 & 0.212 & 0.220 & 2.855 & 0.280 & 0.135 & 0.194 \\
  1 & 3.446 & 1.016 & 0.198 & 0.205 & 3.257 & 0.348 & 0.030 & 0.034 \\
    3 & 3.650 & 1.093 & 0.185 & 0.214 & 3.518 & 0.376 & 0.025 & 0.035 \\
    10 & 3.849 & 1.108 & 0.186 & 0.215 & 3.705 & 0.380 & 0.024 & 0.024 \\
    
 \hline\hline
\kappa & \multicolumn{4}{c|}{m=2^{4}} & \multicolumn{4}{c|}{m=2^{6}} \\ \hline
1/2 & 3.248 & 1.693 & 0.800 & 0.801 & 8.102 & 0.878 & 0.562 & 0.744 \\
  1  & 3.618 & 2.030 & 0.747 & 0.769 & 9.190 & 1.029 & 0.214 & 0.218 \\
  3  & 3.822 & 2.201 & 0.646 & 0.753 & 9.814 & 1.125 & 0.166 & 0.181 \\
  10  & 3.941 & 2.274 & 0.729 & 0.844 & 10.016 & 1.153 & 0.170 & 0.181 \\

 \hline\hline
\kappa & \multicolumn{4}{c|}{m=2^{6}} & \multicolumn{4}{c|}{m=2^{9}} \\ \hline
1/2  & 9.061 & 3.504 & 3.119 & 2.775 & 16.961 & 1.224 & 1.613 & 1.510 \\
  1  & 9.784 & 4.194 & 2.943 & 2.850 & 18.686 & 1.548 & 1.758 & 1.785 \\
  3  & 9.751 & 4.660 & 2.085 & 2.368 & 18.878 & 2.220 & 1.401 & 1.568 \\
  10  & 9.772 & 4.861 & 2.276 & 2.546 & 18.739 & 2.700 & 0.992 & 1.126 \\
 
 \hline
\end{tabu}$$
}
\end{table}

\begin{table}
\caption{The scaled MISE values ($\times 100$) and sd values ($\times 100$) for the estimators}
{\fontsize{10pt}{9pt}\selectfont
$$\begin{tabu}[c]{|c||cccc|cccc|}
   \hline
    {\rm inv. Burr} & \multicolumn{4}{c|}{n=2^{8}} & \multicolumn{4}{c|}{n=2^{12}} \\
  & {\rm PE} & {\rm sd} & {\rm NE} & {\rm sd} & {\rm PE} & {\rm sd} & {\rm NE} & {\rm sd} \\  \hline\hline 
    
 c,\ell & \multicolumn{4}{c|}{m=2^{2}} & \multicolumn{4}{c|}{m=2^{3}}  \\ \hline
3, 2 & 0.220 & 0.201 & 0.188 & 0.213 & 0.056 & 0.037 & 0.023 & 0.046 \\
1,2 & 0.118 & 0.145 & 0.145 & 0.429 & 0.007 & 0.009 & 0.007 & 0.021 \\
1/2, 2 & 11.670 & 10.232 & 3.967 & 0.178 & 10.775 & 1.784 & 8.578 & 0.006 \\
3,1 & 0.186 & 0.198 & 0.178 & 0.274 & 0.034 & 0.067 & 0.023 & 0.027 \\
1,1 & 1.393 & 1.761 & 3.800 & 2.795 & 1.188 & 1.247 & 9.573 & 4.699 \\
1/2,1 & 14.297 & 7.863 & 2.112 & 0.172 & 16.025 & 0.958 & 4.562 & 0.008 \\ 
3,1/3 & 0.329 & 0.531 & 0.152 & 0.136 & 0.170 & 0.441 & 0.015 & 0.015 \\
1,1/3 & 2.322 & 2.785 & 2.581 & 1.570 & 5.398 & 3.153 & 5.498 & 2.307 \\
1/2,1/3 & 19.177 & 12.965 & 1.080 & 0.091 & 25.724 & 13.055 & 2.576 & 0.005 \\

 \hline\hline
 c,\ell & \multicolumn{4}{c|}{m=2^{4}} & \multicolumn{4}{c|}{m=2^{6}}  \\ \hline

3,2 & 0.805 & 0.827 & 1.047 & 1.419 & 0.174 & 0.196 & 0.224 & 0.768 \\
1,2 & 0.802 & 0.722 & 20.412 & 4.993 & 0.165 & 0.169 & 27.016 & 0.823 \\
1/2,2 & 3.151 & 2.883 & 10.800 & 0.000 & 1.049 & 1.469 & 13.800 & 0.000 \\
3,1 & 1.005 & 1.006 & 1.030 & 1.982 & 0.168 & 0.168 & 0.345 & 0.864 \\
1,1 & 1.171 & 1.493 & 15.357 & 3.087 & 0.457 & 0.813 & 18.382 & 0.064 \\
1/2,1 & 9.562 & 10.879 & 5.020 & 0.001 & 10.410 & 10.025 & 5.260 & 0.000 \\
3,1/3 & 1.345 & 1.680 & 1.635 & 1.398 & 0.499 & 0.838 & 1.758 & 1.006 \\
 1,1/3 & 2.340 & 2.747 & 8.806 & 0.720 & 4.089 & 5.180 & 9.120 & 0.001 \\
1/2,1/3 & 22.003 & 21.018 & 2.630 & 0.000 & 21.622 & 12.551 & 2.560 & 0.000 \\ 

 \hline\hline
 c,\ell & \multicolumn{4}{c|}{m=2^{6}} & \multicolumn{4}{c|}{m=2^{9}}  \\ \hline

3,2 & 5.170 & 3.688 & 4.534 & 5.414 & 2.215 & 2.176 & 1.273 & 2.898 \\
1,2 & 4.964 & 5.610 & 27.141 & 0.612 & 1.861 & 1.704 & 29.300 & 0.000 \\
1/2,2 & 3.955 & 4.655 & 13.800 & 0.000 & 2.903 & 3.199 & 16.300 & 0.000 \\
3,1 & 5.394 & 4.581 & 5.413 & 5.139 & 2.268 & 2.519 & 4.159 & 4.650 \\
1,1 & 4.251 & 3.587 & 18.379 & 0.089 & 2.115 & 2.622 & 18.500 & 0.000 \\
1/2,1 & 16.127 & 16.125 & 5.260 & 0.000 & 13.913 & 14.719 & 5.330 & 0.000 \\
3,1/3 & 4.917 & 5.485 & 11.581 & 3.383 & 2.076 & 2.392 & 17.875 & 0.410 \\
1,1/3 & 6.564 & 8.158 & 9.120 & 0.000 & 7.986 & 8.447 & 9.010 & 0.000 \\
1/2,1/3 & 18.610 & 23.975 & 2.560 & 0.000 & 28.425 & 24.743 & 2.540 & 0.000 \\

 \hline
\end{tabu}$$
}
\end{table}

\clearpage

\section{Real data study}

This section considers two real case studies. The first case is the Potomac River peak stream flow (cfs) data for water years (Oct--Sep) 1895--2000 at Point of Rocks, Maryland. The Potomac dataset is available in the \texttt{extRemes} package in R. Figure 1 shows the series of the peak flows ($1 \text{cfs}=0.028317 m^3/s$). The maximal value $480000$ was observed in 1936, and the minimum value was $27800$ in 1969. \cite{smith1987estimating} obtained the fitted GEV to the data (up to 1986) by the maximum likelihood method, where $\widehat{\gamma}=0.42$. \cite{katz2002statistics} argued that the estimated distribution does not necessarily fit the upper tail and obtained fairly strong evidence of a heavy-tail ($\widehat{\gamma}=0.191$ and $p$-value $= 0.002$ for the likelihood ratio test of ${\gamma}=0$) by fitting the GEV to the time series of 106 annual maxima.

We chose the series of $n=100$ peaks from 1901 to 2000. Maximizing the likelihood of the annual peak flows, we obtained $(\widehat{\gamma}_{1^*},\widehat{a}_{1^*},\widehat{b}_{1^*})=(0.200,60000,97100)$, where $1^*$ means that the parameters were estimated by the maximum over each of the years. In a similar manner, $(\widehat{\gamma}_{5^*},\widehat{a}_{5^*},\widehat{b}_{5^*})=(0.847,80300,183000)$, $(\widehat{\gamma}_{10^*},\widehat{a}_{10^*},\widehat{b}_{10^*})=(-0.128,91600,236000)$ and $(\widehat{\gamma}_{20^*},\widehat{a}_{20^*},\widehat{b}_{20^*})=(-0.301,85900,310000)$, where the parameters correspond to the maximum over 5 years, 1 decade and 2 decades, respectively. The occurrence probabilities of the maximum peak flow taking more than some values are estimated, which are given in Table 6. The standard deviations (sd) were estimated by the leave-one-out cross-validation. The results of the extreme value index estimation imply light-tailness. Comparing the estimated probabilities of the peak flows, we see that NE seems optimistic especially for 16000--240000 flows. NE for higher risks returns the high probabilities as the forecast period increases. 

\begin{figure}
\caption{The annual peak flow of the Potomac River at Point of Rocks, MD, USA, 1895--2000}
\begin{center}
\includegraphics[height=7cm]{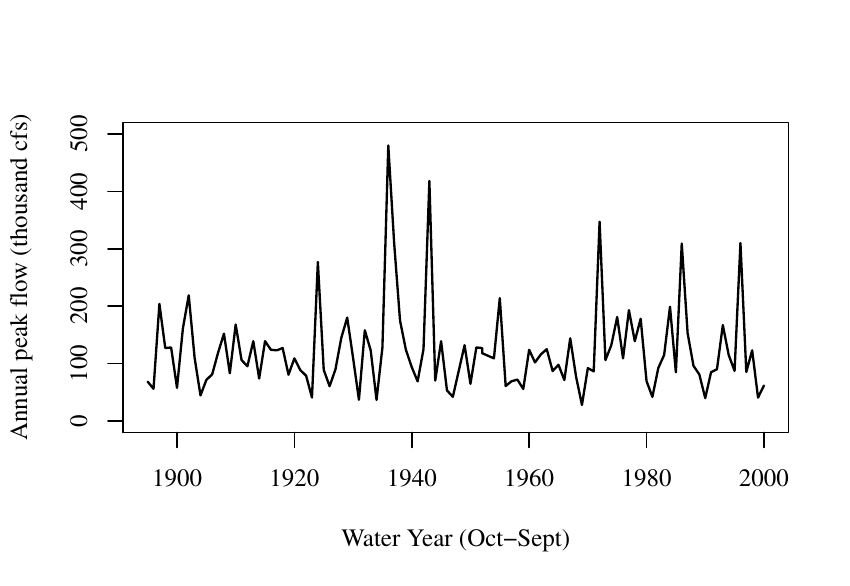}
\end{center}
\end{figure}

\begin{figure}
\begin{center}
\caption{The losses over one million Danish Krone (DKK)} 
\includegraphics[height=7cm]{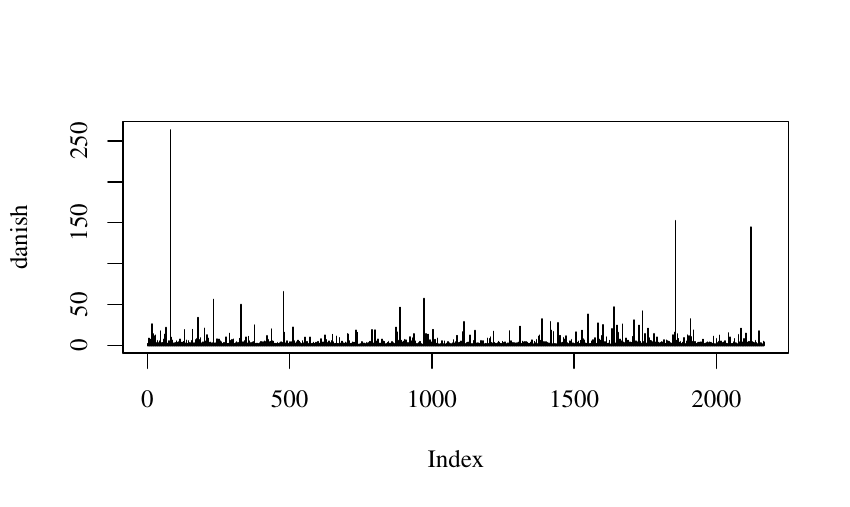}
\end{center}
\end{figure}

\begin{table}
\caption{Estimated probabilities of peak flows of the Potomac River}
$$\begin{tabu}[c]{|c|ccccc|}
  \hline
& \multicolumn{5}{c|}{\rm annual ~ Peak ~  flows} \\
 & 160000 & 240000 & 320000 & 400000 & 480000 \\ \hline

 {\rm PE} & 0.320 & 0.133 & 0.060 & 0.030 & 0.016 \\ 
 {\rm sd} & 0.051 & 0.029 & 0.015 & 0.009 & 0.005 \\
 {\rm NE} & 0.186 & 0.072 & 0.038 & 0.018 & 0.005 \\ 
 {\rm sd} & 0.009 & 0.003 & 0.002 & 0.001 & 0.000 \\
   \hline
& \multicolumn{5}{c|}{\rm 5 ~ years ~ Peak ~  flows} \\
 & 160000 & 240000 & 320000 & 400000 & 480000 \\ \hline

 {\rm PE} & 0.750 & 0.436 & 0.294 & 0.217 & 0.171 \\ 
  {\rm sd} & 0.051 & 0.084 & 0.072 & 0.060 & 0.051 \\
 {\rm NE} & 0.643 & 0.311 & 0.175 & 0.089 & 0.025 \\ 
  {\rm sd} & 0.016 & 0.010 & 0.009 & 0.004 & 0.000 \\
   \hline
& \multicolumn{5}{c|}{\rm 1 ~ decade ~ Peak ~ flows} \\
 & 160000 & 240000 & 320000 & 400000 & 480000 \\ \hline

 {\rm PE} & 0.888 & 0.615 & 0.313 & 0.122 & 0.038 \\ 
   {\rm sd} & 0.013 & 0.021 & 0.013 & 0.004 & 0.005 \\
 {\rm NE} & 0.873 & 0.526 & 0.320 & 0.170 & 0.049 \\ 
  {\rm sd} & 0.010 & 0.013 & 0.015 & 0.008 & 0.000 \\
    \hline
& \multicolumn{5}{c|}{\rm 2 ~ decades ~ Peak ~ flows} \\
 & 160000 & 240000 & 320000 & 400000 & 480000 \\ \hline

 {\rm PE} & 0.983 & 0.873 & 0.587 & 0.246 & 0.048 \\ 
    {\rm sd} & 0.007 & 0.039 & 0.025 & 0.014 & 0.030 \\
 {\rm NE} & 0.984 & 0.775 & 0.537 & 0.311 & 0.095 \\ 
 {\rm sd} & 0.002 & 0.013 & 0.020 & 0.014 & 0.000 \\
 \hline
\end{tabu}$$
\end{table}

\begin{table}
\caption{Estimated probabilities of the fire insurance loss}
$$\begin{tabu}[c]{|c|ccccc|}
  \hline
& \multicolumn{5}{c|}{\rm 1 ~ loss} \\
 & ~~x_{n,5}~~ & ~~x_{n,4}~~ & ~~x_{n,3}~~ & ~~x_{n,2}~~ & ~~x_{n,1}~~ \\ \hline

 {\rm PE} & 0.005 & 0.003 & 0.002 & 0.002 & 0.001 \\ 
     {\rm sd} & 0.000 & 0.000 & 0.000 & 0.000 & 0.000 \\
 {\rm NE} & 0.001 & 0.001 & 0.000 & 0.000 & 0.000 \\ 
     {\rm sd} & 0.000 & 0.000 & 0.000 & 0.000 & 0.000 \\ 
   \hline
& \multicolumn{5}{c|}{\rm 10 ~ losses} \\
 & x_{n,5} & x_{n,4} & x_{n,3} & x_{n,2} & x_{n,1} \\ \hline
 
  {\rm PE} & 0.018 & 0.010 & 0.007 & 0.005 & 0.004 \\ 
       {\rm sd} & 0.001 & 0.000 & 0.000 & 0.000 & 0.000 \\
 {\rm NE} & 0.014 & 0.014 & 0.005 & 0.005 & 0.002 \\ 
      {\rm sd} & 0.000 & 0.000 & 0.000 & 0.000 & 0.000 \\
    \hline
& \multicolumn{5}{c|}{\rm 30 ~ losses} \\
 & x_{n,5} & x_{n,4} & x_{n,3} & x_{n,2} & x_{n,1} \\ \hline
 
   {\rm PE} & 0.042 & 0.021 & 0.013 & 0.009 & 0.006 \\ 
        {\rm sd} & 0.001 & 0.001 & 0.001 & 0.001 & 0.000 \\
 {\rm NE} & 0.042 & 0.042 & 0.014 & 0.014 & 0.007 \\ 
      {\rm sd} & 0.001 & 0.001 & 0.000 & 0.000 & 0.000 \\
    \hline
& \multicolumn{5}{c|}{\rm 100 ~ losses} \\
 & x_{n,5} & x_{n,4} & x_{n,3} & x_{n,2} & x_{n,1} \\ \hline
 
  {\rm PE} & 0.151 & 0.087 & 0.058 & 0.043 & 0.033 \\ 
       {\rm sd} & 0.002 & 0.001 & 0.001 & 0.001 & 0.001 \\
 {\rm NE} & 0.133 & 0.133 & 0.047 & 0.047 & 0.024 \\ 
      {\rm sd} & 0.002 & 0.002 & 0.001 & 0.001 & 0.000 \\
     \hline
& \multicolumn{5}{c|}{\rm 200 ~ losses} \\
 & x_{n,5} & x_{n,4} & x_{n,3} & x_{n,2} & x_{n,1} \\ \hline
 
  {\rm PE} & 0.256 & 0.150 & 0.101 & 0.074 & 0.057 \\ 
       {\rm sd} & 0.005 & 0.002 & 0.001 & 0.001 & 0.001 \\
 {\rm NE} & 0.249 & 0.249 & 0.091 & 0.091 & 0.047 \\ 
     {\rm sd} & 0.004 & 0.003 & 0.002 & 0.002 & 0.000 \\
 \hline
\end{tabu}$$
\end{table}

Next, we consider Danish Fire Insurance data, which are available in the \texttt{evir} package (see also \citealt{beirlant1996practical}). The data were supplied by Mette Rytgaard of Copenhagen Re. It consists of losses over the one million Danish Krone (DKK) (see \citealt{resnick1997discussion}) collected from 1980 to 1990 inclusive. The values were adjusted for inflation to the 1985 values. Figure 2 shows the series of fire insurance losses, where the maximum loss is 263.25 and the minimum is 1. \cite{mcneil1997estimating} fitted the generalized Pareto distribution to the data and discussed the tail. \cite{resnick1997discussion} reported that there is very little doubt about independence, which supported McNeil's (1997) conclusion. The heavy-tail was detected by \cite{resnick1997discussion} (see also \cite{resnick2007heavy} for further studies). 

We chose the latest $n=2100$ losses and obtained $(\widehat{\gamma}_{1^*},\widehat{a}_{1^*},\widehat{b}_{1^*})=(0.922,0.578,1.47)$, $(\widehat{\gamma}_{10^*},\widehat{a}_{10^*},\widehat{b}_{10^*})=(0.695,3.72,5.87)$, $(\widehat{\gamma}_{30^*},\widehat{a}_{30^*},\widehat{b}_{30^*})=(0.560,8.65,12.7)$, $(\widehat{\gamma}_{100^*},\widehat{a}_{100^*},\widehat{b}_{100^*})=(0.718,16.3,27.3)$ and $(\widehat{\gamma}_{200^*},\widehat{a}_{200^*},\widehat{b}_{200^*})=(0.706,24.7,40.1)$,which is estimated by the 2100 losses, the 210 maxima, 70 maxima, 21 maxima and 11 maxima, respectively. The occurrence probabilities of the maximum loss taking more than $x_{n,1}:=263.25$, $x_{n,2}:= (5/6)x_{n,1} = 219.38$, $x_{n,3}:= (2/3)x_{n,1} = 175.5$, $x_{n,4}:= (1/2)x_{n,1} = 131.63$, and $x_{n,5}:= (1/3)x_{n,1} = 87.75$ were obtained and are summarized in Table 7. The sd values were estimated by the leave-one-out cross-validation. 

Since all the estimated shape parameters are not particularly close to zero, we may not need to employ NE. Comparing the estimated probabilities, we find that the estimated risks are not one-sided. NE tends to estimate the risks at $x_{n,4}$ and $x_{n,2}$ to be higher and return the  higher probabilities for high risks as $m$ increases. 

\section{Conclusion and Discussion}

We have discussed the probability distribution estimation of sample maximum sampling from the three classes of distributions. We have confirmed that there exist the two approaches PE and NE based on the extreme value theory and the kernel smoothing, respectively. This study shows that both of the convergence rates heavily depend on $m$ and its extreme value index $\gamma$. PE is effective for large $m$ s.t. the approximation error $F^m - G_{\bm{\gamma}}$ is small enough, while NE is considered to be effective for small $m$. The performances in finite-sample cases were also demonstrated to heavily depend on not only $m$ but also $\gamma$. It was suggested that NE numerically performs well, especially for $\gamma$ close to zero. 

The accuracy of NE depends on the method of choosing bandwidth $h$. Since the sd values of NE in the simulation study were almost moderate except for the bounded cases, we see that the \cite{altman1995bandwidth}'s bandwidth estimator controlled MISE to some degree. As seen in Section 2, the bandwidth estimators of $\widehat{F}$ are the candidates for the bandwidth $h$ of NE. Comparative investigations of the existing bandwidth estimators of $\widehat{F}$ are needed; however, we note that the condition of the optimal bandwidth of NE for the convergence breaks in many cases. In such a case any bandwidth satisfying $h \to 0$ can be used under some conditions (see Remark 4). We need to study the bandwidth selection from this point of view. 

A way to improve the convergence rate of NE is to reduce the asymptotic bias. Since the asymptotic bias of NE coincides with $m$ times that of $\widehat{F}$, the existing methods for bias reduction of the kernel distribution estimator can be used (see, e.g. \citealt{jones1993simple}). Especially, NE for the bounded class of distributions with $\lim_{n\to \infty} h(x^* -x)^{-1} >0$ has the so-called boundary bias , which should be reduced by the recently developing approach e.g. \cite{kolacek2011generalized} and \cite{tenreiro2013boundary}.

It is possible that there exist some nonparametric approaches to SMD estimation. Examples include the `naive' kernel distribution estimator $\bar{F}_{(m)}$ using the block maxima $Y_j := \max\{X_{m(j-1)+1},X_{m(j-1)+2},\cdots,X_{mj}\} (j=1,\cdots,n/m)$, whose distribution is SMD itself $F^m$. Corollary 3 states that NE with the optimal bandwidth has the asymptotic bias of order $(m/n)^{2/3}$ and the variance of order $(m/n)$. The orders are considered to be same as those of $\bar{F}_{(m)}$, since the asymptotic bias of the optimal naive kernel distribution estimator with sample size $N:=(n/m)$ is of order $N^{-2/3}$ and the variance is of order $N^{-1}$. To study the difference between asymptotic properties of the nonparametric approaches may be an important future work. 

The simulation study reveals NE works poorly for $\gamma$ far from zero. There are existing studies on constructing nonparametric density estimators for heavy-tailed data: \cite{wand1991transformations}, \cite{bolance2003kernel}, \cite{maiboroda2004estimation}, \cite{buch2005kernel} and so on. To devise a nonparametric SMD estimator for heavy-tailed data merits further research. 

NE is expected to have the consistency even if $F$ does not belong to any maximum domain of attraction. For the cases, we cannot apply the fitting approach PE. Both theoretical and numerical studies on the cases are desired. 

Considering the above, each of PE and NE has its own strengths. One of the next important issues is devising semiparametric estimators that have both of the properties of PE and NE. 

\cite{grinevigh1994domains} introduced the class of max-semistable laws, which hold even if the underlying distribution $F$ is discrete or multimodal. The usual max-stable laws do not hold in the cases, which means PE cannot be applied. \cite{temido2003max} and \cite{canto2011looking} studied more details (see also \cite{pancheva2010max}). Extending this study to wider classes is a topic for future work.

\bibliographystyle{natbib}
\bibliography{ref}

\section*{Appendices: on fitted estimation of SMD based on extreme value theory}

The limiting distributions of the Hall class, the Weibull class and the bounded class of distributions are the Fr\'{e}chet, Gumbel and Weibull type of GEV, respectively (see \cite{beirlant2004statistics}) with
\begin{align*}
\gamma :=&
\begin{dcases}
\alpha^{-1} ~~~ &\text{for the Hall class} \\
0 ~~~ &\text{for the Weibull class} \\
\mu^{-1} ~~~ &\text{for the bounded class},
\end{dcases} \\
a_m :=&
\begin{dcases}
\gamma (Am)^{\gamma} ~~~ &\text{for the Hall class} \\
\kappa^{-1} C^{-1/\kappa} (\ln m)^{-\theta} ~~~ &\text{for the Weibull class} \\
-\gamma (D m)^{\gamma} ~~~ &\text{for the bounded class},
\end{dcases} \\
\intertext{$\theta :=1-(1/\kappa)$, and} \\
b_m :=&
\begin{dcases}
(Am)^{\gamma} ~~~ &\text{for the Hall class} \\
(C^{-1} {\ln m})^{1/\kappa} ~~~ &\text{for the Weibull class} \\
x^* - (D m)^{\gamma} ~~~ &\text{for the bounded class}.
\end{dcases}
\end{align*}

This study considers the MLE based on the block maxima method (\cite{dombry2019maximum}). We suppose that the block size is $k := k_n (\to \infty)$ and $\exists N \in \mathbb{N}$ s.t. $n=N\times k$. In the Appendices the block size is explicitly expressed in the subscript like as $\bm{\gamma}_k$. Set
$$\widehat{\bm{\gamma}}_k := \mathop{\rm arg~max}\limits_{\bm{\gamma}_k} \sum_{j=1}^N \ln g_{\bm{\gamma}_k}(Y_j),$$
where $Y_j := \max\{X_{k(j-1)+1},X_{k(j-1)+2},\cdots,X_{kj}\} (j=1,\cdots,N)$ and $g_{\bm{\gamma}}$ is the density function of $G_{\bm{\gamma}}$. We can set any value as $k$ satisfying the following assumption, that means $k=m$ is not necessary. 

\begin{assumption}
Any one of the following {\rm (i)}--{\rm (iii)} holds
\begin{align*}
&{\rm (i)} ~~~ (M_n \vee K_n) \to 0\\
&{\rm (ii)} ~~~ (M_n \wedge K_n) \to \infty \\
&{\rm (iii)} ~~~ \exists\delta>0 ~~ {\rm s.t.} ~~  M_n \to \delta ~~ {\rm and} ~~ K_n \to \delta,
\end{align*}
where 
\begin{align*}
K_n :=
\begin{dcases}
A k x^{-\alpha} ~~~ &\text{for the Hall class} \\
k^{\kappa} \exp (-\kappa C^{1/\kappa} (\ln k)^{\theta} x ) ~~~ &\text{for the Weibull class} \\
Dk (x^* -x)^{-\mu} ~~~ &\text{for the bounded class}.
\end{dcases}
\end{align*}
\end{assumption}

Assumption 1 ensures the convergence of the approximation of SMD to GEV (see Proposition 1 below). 

\begin{proposition}
Under Assumption 1
$$\tau_n := F^m(x) - G_{{\bm{\gamma}}_k}(x) \to 0.$$
\end{proposition}

As noted later, $k$ controls the bias-variance trade of  the estimator $\widehat{\bm{\gamma}}_k$. After \cite{prescott1980maximum}, \cite{dombry2015existence} proved the existence and consistency of the MLE in the block maxima framework under a practical condition. After \cite{bucher2017maximum}, \cite{dombry2019maximum} provided the conditions that the following scaled version of the MLE
$$\widehat{\bm{\gamma}}_k^*=\begin{pmatrix}
1 & & \\
 & a_k^{-1} & \\
 & & a_k^{-1}
\end{pmatrix}
\widehat{\bm{\gamma}}_k$$ 
exists and is asymptotically normal with a nontrivial bias, $\sqrt{N}$-consistent and efficient under the following Assumption 2. 

\begin{assumption}
$\exists \lambda \in \mathbb{R}$ s.t. $\lambda_n \to \lambda$, where 
$$\lambda_n := \sqrt{k} \times
\begin{dcases}
m^{-\beta} ~~~ &\text{for the Hall class}, \\
(\ln m)^{-1} ~~~ &\text{for the Weibull class}, \\
m^{\sigma} ~~~ &\text{for the bounded class}.
\end{dcases}
$$
\end{assumption}
Assumption 2 ensures the convergence of the asymptotic bias of the MLE $\widehat{\bm{\gamma}}_k$. Note that $k \sim m$ for the Weibull class means $\lambda_n \to \infty$ i.e. Assumption 2 broken. 

We have the following proposition on the accuracy of the pointwise estimation of SMD, which is a direct consequence of Theorem 2.2 in \cite{dombry2019maximum}.

\begin{proposition}
Under Assumptions 1--2,
$$\mathbb{E}[(F^m(x) - G_{\widehat{\bm{\gamma}}_k}(x))^2] \sim (\tau_n - N^{-1/2} \lambda_n \bm{\eta}_n^{\mathsf{T}} I_{0}^{-1} \bm{b})^2 + N^{-1} (\bm{\eta}_n^{\mathsf{T}} I_{0}^{-1} \bm{\eta}_n),$$ 
where $\bm{b}$ and the Fisher information matrix $I_{0}$ are given in \cite{dombry2019maximum}, and 
\begin{align*}
\bm{\eta}_n := -\exp(-K_n) K_n (1- K_n^{\gamma} + \gamma \ln K_n, K_n^{\gamma} \left(\frac{K_n^{-\gamma} -1}{\gamma}\right), K_n^{\gamma})^{\mathsf{T}}.
\end{align*}
Furthermore, 
$$N^{1/2} (\bm{\eta}_n^{\mathsf{T}} I_{0}^{-1} \bm{\eta}_n)^{-1/2} \{F^m(x) -\widehat{F}^m(x) - (\tau_n - N^{-1/2} \lambda_n \bm{\eta}_n^{\mathsf{T}} I_{0}^{-1} \bm{b})\}$$
converges in distribution to the standard normal distribution if the MSE converges to zero.
\end{proposition}

The three bias terms $F^m(x) -G_{\bm{\gamma}_m}(x)$, $G_{\bm{\gamma}_m}(x) - G_{\bm{\gamma}_k}(x)$ and $\lambda_n$ appears, as seen from the following decomposition:
$$F^m(x) - G_{\widehat{\bm{\gamma}}_k}(x) = \{F^m(x) -G_{\bm{\gamma}_m}(x)\} -\{G_{\bm{\gamma}_m}(x) - G_{\bm{\gamma}_k}(x)\} - \{G_{\bm{\gamma}_k}(x) - G_{\widehat{\bm{\gamma}}_k}(x)\}.$$
$\{G_{\bm{\gamma}_m}(x) -G_{\bm{\gamma}_k}(x)\}$ requires the asymptotic equivalence of $k$ and $m$ for its convergence. The natural choice $k \sim m$ requires the Hall class or the bounded class. Although a large $k$ yields a large variance of $\widehat{\bm{\gamma}}_k$, the precise approximation of $G_{\bm{\gamma}_m}(x)$ to $F^m(x)$ (i.e. $\tau_n$) requires a large $m$. In summary, a small $m$ or an excessively large $m$ compared to $n$ makes the convergence rate slow. $\lambda_n$ is the asymptotic bias of estimator $\widehat{\bm{\gamma}}_k$, i.e., $G_{\widehat{\bm{\gamma}}_k}(x)$. Numerical examples of the convergence rate are given in Table 1 in Section 3.

\section*{Appendices: Proofs}

\noindent\textbf{Proof of Proposition 1 for the Hall class} 

It follows from $a_k = \alpha^{-1} (Ak)^{\gamma}$ and $b_k = (Ak)^{\gamma}$ that
$$G_{\bm{\gamma}_k}(x) = \exp \left(- \left(1+ \gamma \frac{x- b_k}{a_k} \right)^{-\alpha}\right) = \exp\left(- A k x^{-\alpha} \right).$$
It follows from
$$\ln F^m(x) \sim -M_{n} \left\{ 1 + B x^{-\beta} + \frac{1}{2} A x^{-\alpha} \right\}$$
that
$$F^m(x) \sim \exp(-M_{n}) \exp\left[-M_{n} \left\{ B x^{-\beta} + \frac{1}{2} A x^{-\alpha} \right\}\right],$$
where $M_{n}=A m x^{-\alpha}$ and $K_{n}=A k x^{-\alpha}$. Then, 
$$\tau_n \sim \exp\left(- M_{n} \right)\exp\left[-M_{n} \left\{ B x^{-\beta} + \frac{1}{2} A x^{-\alpha} \right\}\right] -\exp\left(- K_{n} \right).$$
It is immediately seen that $\tau_n \to 0$ if $(M_n \vee K_n) \to 0$ or $(M_n \wedge K_n) \to \infty$. When there exists some positive constant $\delta$ s.t. $M_n - \delta =: \epsilon_{n} \to 0$, $\tau_n$ converges if and only if $K_n - \delta =: \widetilde{\epsilon_{n}} \to 0$. Then,
$$\tau_n \sim \exp(- \delta) \biggm\{ \exp \left[- \epsilon_{n} -\delta \left\{ B m^{-\beta \gamma} + \frac{1}{2} A m^{-1} \right\} \right] -\exp\left(- \widetilde{\epsilon}_{n} \right) \biggm\}$$
converging with the rate $(m^{-\beta \gamma} \vee m^{-1} \vee {\epsilon}_{n} \vee \widetilde{\epsilon}_{n})$. Proposition 1 for the Hall class has been proved. \\

\vspace{\baselineskip}
\noindent\textbf{Proof of Proposition 1 for the Weibull class} 

It follows from $a_k = \kappa^{-1} C^{-1/\kappa} (\ln k)^{-\theta}$ and $b_k = (C^{-1} {\ln k})^{1/\kappa}$ that
\begin{align*}
G_{\bm{\gamma}_k}(x) =& \exp \left(- \exp \left( -\frac{x -b_{k}}{a_{k}} \right)\right) = \exp [ -k^{\kappa} \exp (-\kappa C^{1/\kappa} (\ln k)^{\theta} x)].
\end{align*}
It holds that
\begin{align*}
F^m(x) =& \left(1 -\exp\{-C x^{\kappa}\} \right)^m \sim \exp\left[ m \ln \left(1 -\exp\{-C x^{\kappa}\} \right) \right] \\
\sim& \exp\left( -m \exp\{-C x^{\kappa}\} - 2^{-1}m \exp\{-2C x^{\kappa}\} \right) \\
\sim&
\begin{dcases}
1 -M_n + 2^{-1} M_n^2 (1 -m^{-1}) ~~~ \text{for} ~~~ M_n \to 0\\
\exp\left( -M_n \right)  ~~~ \text{else},
\end{dcases}
\end{align*}
where $M_n= m \exp(-C x^{\kappa})$. It follows that the condition of $\tau_n \to 0$ is (i) $(M_n \vee K_n) \to 0$, (ii) $(M_n \wedge K_n) \to \infty$ or (iii) $\exists\delta>0$ s.t. $M_n \to \delta$ and $K_n \to \delta$, where $K_n=k^{\kappa} \exp (-\kappa C^{1/\kappa} (\ln k)^{\theta} x )$.

\clearpage

\vspace{\baselineskip}
\noindent\textbf{Proof of Proposition 1 for the bounded class} 

It holds that
\begin{align*}
G_{\bm{\gamma}_k}(x) & = \exp \left(- \left(1+ \gamma \frac{x- b_k}{a_k} \right)^{-\mu}\right) = \exp\left( -Dk (x^* -x)^{-\mu}\right) \\
\ln F^m(x) &= -m \{(1-F(x)) +\frac{1}{2} (1-F(x))^2 + O((1-F(x))^3)\} \\
F^m(x) &= \exp(-m \{(1-F(x)) +\frac{1}{2} (1-F(x))^2 + O((1-F(x))^3)\}).
\end{align*}
When $(M_n \vee K_n) =o(1)$ and $(M_n \wedge K_n) \to \infty$, 
\begin{align*}
{\tau_n} =& K_n -M_n +o(M_n \vee K_n) \to 0, \\
{\tau_n} =& \exp\left(- M_{n} + o(M_n)\right) - \exp\left(- K_{n} +o(K_n)\right) \to 0
\end{align*}
respectively. When there exists some positive constant $\delta$ s.t. $M_n - \delta =: \epsilon_{n} \to 0$, ${\tau}_n$ converges if and only if $K_n - \delta =: {\epsilon_{n}} \to 0$. Then,
\begin{align*}
&{\tau}_n - \exp(-\delta) \left[\exp\left(-\epsilon_n -\delta\left\{D^{-1} E (x^* -x)^{-\mu\sigma} + \frac{\delta}{2}m^{-1} \right\} \right) -\exp(-{\epsilon}_n)\right]
\end{align*}
is of order $(m^{-\mu\sigma} \vee m^{-1} \vee \epsilon_n \vee {\epsilon}_n)$. \\

\vspace{\baselineskip}
\noindent\textbf{Proof of Proposition 2 for the Hall class or the bounded class} 

First, we decompose the difference as follows:
\begin{align*}
F^m(x) - G_{\widehat{\bm{\gamma}}_k}(x) =& \left[F^m(x) - G_{\bm{\gamma}_k}(x) \right] - \left[G_{\bm{\gamma}_k}(x) - G_{\widehat{\bm{\gamma}}_k}(x)\right] =: \tau_n + \zeta_n ~~~~~~~~~ (\text{say}).
\end{align*}
It holds that
$$\zeta_n=G_{\bm{\gamma}_k}(x) - G_{\widehat{\bm{\gamma}}_k}(x) = -\frac{\partial}{\partial \bm{\gamma}}G_{\bm{\gamma}}(x) \biggm\vert_{\bm{\gamma}=\widetilde{\bm{\gamma}}_k} (\widehat{\bm{\gamma}}_k -\bm{\gamma}_k),$$
where $\widetilde{\bm{\gamma}}_k = (\widetilde{\gamma}_k,\widetilde{a}_k,\widetilde{b}_k)^{\mathsf{T}}$ is between $\widehat{\bm{\gamma}}_k$ and $\bm{\gamma}_k$ with probability $1$. By calculating the derivative, we have
\begin{align*}
\frac{\partial}{\partial \alpha}G_{\bm{\gamma}}(x) \biggm\vert_{\bm{\gamma}= \widetilde{\bm{\gamma}}_k} =& -\exp\left\{-\left(1+\frac{x -\widetilde{b}_k}{\widetilde{\alpha}_k \widetilde{a}_k}\right)^{-\widetilde{\alpha}_k}\right\} \left(1+\frac{x -\widetilde{b}_k}{\widetilde{\alpha}_k \widetilde{a}_k}\right)^{-\widetilde{\alpha}_k} \\
& \left(\frac{x -\widetilde{b}_k}{\widetilde{\alpha}_k \widetilde{a}_k + x -\widetilde{b}_k}-\ln\left(1+\frac{x -\widetilde{b}_k}{\widetilde{\alpha}_k \widetilde{a}_k}\right)\right),
\end{align*}
where $\widetilde{\alpha}_k := \widetilde{\gamma}_k^{-1}$. Set $z_{k,n}:=1+ (\alpha {a}_k)^{-1} (x -{b}_k) = \gamma a_k^{-1} x$. It follows from
$$\frac{x -\widehat{b}_k}{\widehat{a}_k} - \frac{x -{b}_k}{{a}_k} = \frac{x}{a_k} \left(\frac{{a}_k}{\widehat{a}_k} -1\right) - \left( \frac{\widehat{b}_k /a_k}{\widehat{a}_k /a_k} - \frac{b_k /a_k}{1}\right)$$
that $(\widehat{z}_{k,n}/z_{k,n}) \xrightarrow{p} 1$. Then,
$$s_n^{-1} \frac{\partial}{\partial \alpha} G_{\widetilde{\bm{\gamma}}_{k}}(x) \overset{p}{\to} 1,$$ 
where $s_n:= -\exp(-K_n) K_n  \left(1 -K_n^{\gamma} +\gamma \ln K_n \right)$.

In the same manner, if $(\widehat{z}_{k,n}/ z_{k,n})  \xrightarrow{p} 1$, it follows from
\begin{align*}
\frac{\partial}{\partial a_k}G_{\bm{\gamma}}(x) \biggm\vert_{\bm{\gamma}=\widetilde{\bm{\gamma}}_k} =& - \widetilde{a}_k^{-1} \exp\left\{-\left(1+\frac{x -\widetilde{b}_k}{\widetilde{\alpha}_k \widetilde{a}_k}\right)^{-\widetilde{\alpha}_k}\right\} \left(\frac{x -\widetilde{b}_k}{\widetilde{a}_k}\right) \left(1+\frac{x -\widetilde{b}_k}{\widetilde{\alpha}_k \widetilde{a}_k}\right)^{-1-\widetilde{\alpha}_k} 
\end{align*}
and 
\begin{align*}
\frac{\partial}{\partial b_k}G_{\bm{\gamma}}(x) \biggm\vert_{\bm{\gamma}=\widetilde{\bm{\gamma}}_k} =& -\frac{1}{\widetilde{a}_k} \exp\left\{-\left(1+\frac{x -\widetilde{b}_k}{\widetilde{\alpha}_k \widetilde{a}_k}\right)^{-\widetilde{\alpha}_k}\right\} \left(1+\frac{x -\widetilde{b}_k}{\alpha \widetilde{a}_k}\right)^{-1-\alpha}
\end{align*}
that
$$(t_n^*)^{-1}\frac{\partial}{\partial a_k}G_{\bm{\gamma}}(x) \biggm\vert_{\bm{\gamma}=\widetilde{\bm{\gamma}}_k} \overset{p}{\to} 1 ~~~\text{and} ~~~ (u_n^*)^{-1} \frac{\partial}{\partial b_k}G_{\bm{\gamma}}(x) \biggm\vert_{\bm{\gamma}=\widetilde{\bm{\gamma}}_k} \overset{p}{\to} 1,$$
where $t_n^* := -a_k^{-1} \exp(-K_n) K_n^{1+\gamma} \left(\frac{K_{n}^{-\gamma} -1}{\gamma}\right)$ and $u_n^* :=-a_k^{-1} \exp(-K_n) K_n^{1+\gamma}$.

As seen from Dombry and Ferrira (2019),
$$
\sqrt{N}
\begin{pmatrix}
1 & & \\
 & a_k^{-1} & \\
 & & a_k^{-1}
\end{pmatrix}
(\widehat{\bm{\gamma}}_k -\bm{\gamma}_k) := \bm{N}$$
converges in distribution to the normal with the mean $\lambda I_{0}^{-1} \bm{b}$ and variance $I_{0}^{-1}$. Thus, we see $\zeta_n$ is asymptotically equivalent in distribution to $-N^{-1/2} \bm{\eta}_n^{\mathsf{T}} \bm{N}$, where $\bm{\eta}_n := (s_n, t_n, u_n)^{\mathsf{T}}$, $t_n:=a_k t_n^*$ and $u_n:=a_k u_n^*$. Combining the results, Proposition 2 for the Hall class has been proved. Proposition 2 for the bounded class is proved in the same manner.\\

\vspace{\baselineskip}
\noindent\textbf{Proof of Proposition 2 for the Weibull class} 

$$\zeta_n :=G_{\bm{\gamma}_k}(x) - G_{\widehat{\bm{\gamma}}_k}(x)= -\frac{\partial}{\partial \bm{\gamma}}G_{\bm{\gamma}}(x) \biggm\vert_{\bm{\gamma}=\widetilde{\bm{\gamma}}_k} (\widehat{\bm{\gamma}}_k -\bm{\gamma}_k)$$
holds, where $\widetilde{\bm{\gamma}}_k$ is between $\widehat{\bm{\gamma}}_k$ and $\bm{\gamma}_k$ with probability $1$. If $(\widehat{z}_{k,n}/ z_{k,n}) \xrightarrow{p} 1$, we have
$$\frac{\partial}{\partial \alpha}G_{\bm{\gamma}}(x) \biggm\vert_{\bm{\gamma}={\widetilde{\bm{\gamma}}_k}} \overset{p}{\to} 0.$$
It holds that
$$\frac{\partial}{\partial a_k}G_{\bm{\gamma}}(x) \biggm\vert_{\bm{\gamma}={\widetilde{\bm{\gamma}}_k}} = -x_{k,n}^{*} {a}_k^{-2} k^{-\phi_{0} {a}_k^{1/(\kappa -1)} x_{k,n}^{*}} \exp\left[-k^{-\phi_{0} {a}_k^{1/(\kappa -1)} x_{k,n}^{*}}\right] =: t_n^*,$$
where $x_{k,n}^{*} :=x-\{C^{-1} \ln k\}^{\frac{1}{\kappa}}$ and $\phi_{0} :=\kappa^{1/{\kappa}}C^{1/(\kappa -1)}$. It follows from
\begin{align*}
k^{-\phi_{0} {a}_k^{1/(\kappa -1)} x_{k,n}^{*}} =& k^{\kappa} \times k^{-\kappa C^{1/\kappa} (\ln k)^{-1/\kappa} x} \\
=& k^{\kappa} \exp\left[\ln (k^{-\kappa C^{1/\kappa} (\ln k)^{-1/\kappa} x}) \right] =: K_n
\end{align*}
that
$$t_n^* - a_k^{-1} \exp(-K_n) K_n \ln K_n \to 0$$
since $a_k^{-1} x_{k,n}^{*} = -\ln K_n$. Therefore, we see if $(\widehat{z}_{k,n}/ z_{k,n}) \xrightarrow{p} 1$
$$(t_n^*)^{-1} \frac{\partial}{\partial a_k}G_{\bm{\gamma}}(x) \biggm\vert_{\bm{\gamma}=\widetilde{\bm{\gamma}}_k} \overset{p}{\to} 1$$
and
$$(u_n^*)^{-1} \frac{\partial}{\partial b_k}G_{\bm{\gamma}}(x) \biggm\vert_{\bm{\gamma}={\widetilde{\bm{\gamma}}_k}} \overset{p}{\to} 1,$$
where $u_n^* := -{a}_k^{-1} K_n \exp(-K_n)$.

In the same manner as the Proof for the Hall class, we see $\zeta_n$ is asymptotically equivalent in distribution to $-N^{-1/2} \bm{\eta}_n^{\mathsf{T}} \bm{N}$, where $\bm{\eta}_n$ $:=$ $(s_n, t_n, u_n)^{\mathsf{T}}$. $s_n \equiv 0$, $t_n:= a_k t_n^*$ and $u_n:= a_k u_n^*$. Proposition 2 for the Weibull class has been proved. \\

\vspace{\baselineskip}
\noindent\textbf{Proof of Corollary 1} 

Corollary 1 immediately follows from Proposition 2 for $k=m$.

\vspace{\baselineskip}
\noindent\textbf{Proof of Theorem 1 for the Hall class} 

Expanding $F(x -hz)$ we have
\begin{align*}
 & 1 -A (x -hz)^{-\alpha} \left(  1 + B (x -hz)^{-\beta} + o(x^{-\beta}) \right)\\
=& 1 -A \left\{x^{-\alpha} -\alpha h x^{-\alpha-1} z + \frac{\alpha(\alpha+1)}{2} h^2 x^{-\alpha-2} z^2 + O(h^3 x^{-\alpha-3})\right\} \\
& ~~~ \left(  1 + B \left\{x^{-\beta} -\beta h x^{-\beta-1} z + \frac{\beta(\beta+1)}{2} h^2 x^{-\beta-2} z^2 + O(h^3 x^{-\beta-3}) \right\} +  o(x^{-\beta}) \right) \\
\text{and} \\
& \mathbb{E}[\widehat{F}(x)] -F(x) \\
=& \int_{-\infty}^{\infty} w\left(z \right) \{F(x -hz) -F(x)\} {\rm d}z \sim -\frac{A\alpha(\alpha+1)}{2} h^2 x^{-\alpha-2} \int_{-\infty}^{\infty} z^2 w\left(z \right) {\rm d}z.
\end{align*}
It holds that
\begin{align*}
n\mathbb{V}[\widehat{F}(x)] =& \mathbb{E}\left[W^2\left(\frac{x-X_i}{h}\right)\right] - \left(\mathbb{E}\left[W\left(\frac{x-X_i}{h}\right)\right]\right)^2 \\
&\sim {A}{x^{-\alpha}} \left\{1 + B {x}^{-\beta} - A {x}^{-\alpha} -2 \alpha h {x^{-1}} \int z W(z) w(z) {\rm d}z\right\}.
\end{align*}

Set $Z_n:= m (\widehat{F}(x) -\mathbb{E}[\widehat{F}(x)])$ and
$$y_n := m (\mathbb{E}[\widehat{F}(x)] - F(x)) \sim -h^{2} m x^{-\alpha-2} \frac{A\alpha(\alpha+1)}{2} \int z^2 w(z) {\rm d}z.$$
It holds that
\begin{align*}
\widehat{F}^m(x) &= \exp \left(m \ln (\widehat{F}(x)) \right) \\
&= \exp\left(m \left\{\widehat{F}(x) -1 + o_P((\widehat{F}(x) -1)) \right\}\right) \\
&= \exp\left( Z_n + y_n + m(F(x) -1) + o_P(m(\widehat{F}(x) -1)) \right) \\
&= \exp\left( Z_n + y_n + o_P(m(\widehat{F}(x) -1)) \right) \exp(-M_{n}),
\end{align*}
where $M_n := Am x^{-\alpha}$ and 
\begin{align*}
\widehat{F}(x) -1 =& (\widehat{F}(x) -\mathbb{E}[\widehat{F}(x)]) + (\mathbb{E}[\widehat{F}(x)] -1) \\
=& O_{P}(n^{-1/2} x^{-\alpha/2}) + (\mathbb{E}[\widehat{F}(x)] -F(x)) + (F(x) -1).
\end{align*}
It follows from $m(\widehat{F}(x) -1)^2 = m x^{-\alpha} O_{P}\left(x^{-\alpha} + n^{-1} \right)$ that
\begin{align*}
F^m(x) - \widehat{F}^m(x) &\sim F^m(x) - \exp\left(-M_n + Z_n + y_n + m x^{-\alpha} O_{P}\left(x^{-\alpha} + n^{-1} \right) \right) ) \\
&\sim F^m(x) (1 - \exp(Z_n +y_n)) \\
& \sim mF^m(x) \{(\widehat{F}(x) -\mathbb{E}[\widehat{F}(x)]) + (\mathbb{E}[\widehat{F}(x)] - F(x)) \},
\end{align*}
where $Z_n = O_{P}\left(m x^{-\alpha/2} n^{-1/2} \right)$. Combining the results, Theorem 1 for the Hall class has been proved.\\

\vspace{\baselineskip}
\noindent\textbf{Proof of Theorem 1 for the Weibull class} 

There exist some $c_0 >0$ such that
$$\left\vert\mathbb{E}[\widehat{F}(x)] -F(x) - h^2 \frac{f'(x)}{2} \int z^2 w(z) {\rm d}z\right\vert < c_0 h^3 \exp(-C x^{\kappa}),$$
where
$$f'(x) = \kappa C x^{\kappa-2} (- \kappa C x^{\kappa} + \kappa -1) \exp(-C x^{\kappa}).$$
For $h x^{\kappa-1} \to 0$, the asymptotic variance is given by 
$$\mathbb{V}[\widehat{F}(x)] \sim n^{-1} \exp(-C x^{\kappa}) \{1 -2 \kappa C h x^{\kappa-1} \int z W(z) w(z) {\rm d}z \}.$$
Following the proof for the Hall class, Theorem 1 for the Weibull class is proved. \\

\vspace{\baselineskip}
\noindent\textbf{Proof of Theorem 1 for the bounded class} 

It follows from
$$F(x -hz) = 1 -(x^* -x +hz)^{-\mu} \left\{D + E (x^* -x +hz)^{-\sigma} + o((x^* -x +hz)^{-\sigma}) \right\}$$
that $F(x -hz)$ asymptotically equals
\begin{align*}
&1 - \left\{(x^* -x)^{-\mu} -\mu (x^* -x)^{-\mu-1} hz +\frac{\mu (\mu+1)}{2} (x^* -x)^{-\mu-2} (hz)^2 \right\} \\
& \Bigm[D + E \left\{(x^* -x)^{-\sigma} + -\sigma (x^* -x)^{-\sigma-1} hz +\frac{-\sigma (-\sigma-1)}{2} (x^* -x)^{-\sigma-2} (hz)^2 \right\} \Bigm]
\end{align*}
for $h (x^* -x)^{-1} \to 0$ or asymptotically equals
\begin{align*}
&1 - \left\{(hz)^{-\mu} -\mu (hz)^{-\mu-1} (x^* -x) +\frac{\mu (\mu+1)}{2} (hz)^{-\mu-2} (x^* -x)^2 \right\} \\
& \Bigm[D + E \left\{(hz)^{-\sigma} +-\sigma (hz)^{-\sigma-1} (x^* -x) +\frac{-\sigma (-\sigma-1)}{2} (hz)^{-\sigma-2} (x^* -x)^2 \right\} \Bigm]
\end{align*}
for $h (x^* -x)^{-1} \to \infty$. $h (x^* -x)^{-1} \to 0$ means ${\rm supp}(w) \subset [h^{-1}(x - x^*), h^{-1}(x - x_*)]$ for large enough $n$, where $x_*$ is the lower endpoint of $f$. $h (x^* -x)^{-1} \to \infty$ means $(x^* -x)^{-\mu}=o(h^{-\mu})$. Therefore, $\mathbb{E}[\widehat{F}(x)] -F(x)$ is asymptotically
$$
\begin{dcases}
-\frac{D \mu (\mu+1)}{2} (x^* -x)^{-\mu-2} h^2 \int z^2 w(z) {\rm d}z~~~ \text{for} ~~~ h (x^* -x)^{-1} \to 0\\
D h^{-\mu} \int_0^{\infty}z^{-\mu}w(z) {\rm d}z +\frac{1}{2}(x^* -x)^{-\mu} ~~~ \text{for} ~~~ h (x^* -x)^{-1} \to \infty,
\end{dcases}
$$
if the support of $w$ is bounded. In a similar manner, $n\mathbb{V}[\widehat{F}(x)]$ is asymptotically
$$
\begin{dcases}
D(x^* -x)^{-\mu} \left(1 -2\mu (x^* -x)^{-1} h \int W(z) w(z) z {\rm d}z\right) ~~~ \text{for} ~~~ h (x^* -x)^{-1} \to 0\\
2 D h^{-\mu} \int_0^{\infty}z^{-\mu} w(z)\{W(z) -1\} {\rm d}z ~~~ \text{for} ~~~ h (x^* -x)^{-1} \to \infty.
\end{dcases}
$$
Applying the proof for the Hall class, Theorem 1 for the bounded class has been proved. \\

\vspace{\baselineskip}
\noindent\textbf{Proof of Corollary 2}

It follows from the proof of Theorem 1 that $\exp(2M_n) \mathbb{E}[(F^m(x) -\widehat{F}^m(x))^2]$ is asymptotically 
$$M_n^{2} h^{4} \frac{\xi_n^2}{4} \left(\int z^2 w(z) {\rm d}z\right)^2 + \frac{m}{n} \left(M_n -2 m h f(x) \int z W(z) w(z) {\rm d}z \right),$$
where 
\begin{align*}
\xi_n:=&\begin{dcases}
\alpha(\alpha+1) x^{-2} ~~~ &\text{for the Hall class}\\
\kappa^2 C^2 x^{2\kappa-2} ~~~ &\text{for the Weibull class},\\
\mu(\mu+1) x^{-2} ~~~ &\text{for the bounded class}\\
\end{dcases} \\
f(x)\sim&\begin{dcases}
A \alpha x^{-\alpha-1} ~~~ &\text{for the Hall class}\\
\kappa C x^{\kappa-1} \exp(-C x^{\kappa}) ~~~ &\text{for the Weibull class},\\
-D \mu x^{-\mu-1} ~~~ &\text{for the bounded class}.
\end{dcases}
\end{align*}
By differencing the MSE with respect to $h$ and summing up, we see that the optimal bandwidth 
$$h= \left( 2 \xi_n^{-2} \psi_n n^{-1} \frac{\int z W(z) w(z) {\rm d}z}{\left(\int z^2 w(z) {\rm d}z\right)^2}\right)^{1/3},$$
where
$$
\psi_n:=\begin{dcases}
A^{-1} \alpha x^{\alpha-1} ~~~ &\text{for the Hall class}\\
\kappa C x^{\kappa-1} \exp(C x^{\kappa}) ~~~ &\text{for the Weibull class},\\
-D^{-1} \mu x^{\mu-1} ~~~ &\text{for the bounded class}.
\end{dcases}
$$

\end{document}